\documentclass[twocolumn]{svjour3}

\smartqed
\usepackage{mathptmx}

\usepackage{amsmath,amsfonts,amssymb}

\usepackage{graphicx}
\graphicspath{{pics/layer-2d/},{pics/random-2d/},{pics/chap4/}}
\usepackage{subcaption}
\newcommand{\avg}[1]{\{\!\!\{#1\}\!\!\}}
\newcommand{\jmp}[1]{[\![#1]\!]}

\newcommand\grad{\nabla}

\renewcommand\div{\nabla\cdot}
\newcommand\Div{\textrm{div}}
\newcommand\Curl{\textrm{curl}}
\newcommand\spn{\textrm{span}}

\newcommand\En{\mathsf E}

\newcommand\Po{\mathbb P}
\newcommand\Qu{\mathbb Q}
 
\renewcommand\Re{\mathbb R}
\newcommand\Su{\mathbb S}

\newcommand\F{{\mathbf F}}
\newcommand\V{{\mathbf V}}
\let\vv=\v
\renewcommand\v{{\mathbf v}}
\newcommand\n{{\pmb\nu}}

\renewcommand\u{{\mathbf u}}
\newcommand\x{{\mathbf x}}

\newcommand\cC{{\cal C}}
\newcommand\cE{{\cal E}}
\newcommand\cF{{\cal F}}
\newcommand\cO{{\cal O}}
\newcommand\cP{{\cal P}}
\newcommand\cT{{\cal T}}

\newcommand\red{\text{\rm red}}
\newcommand\AC{\text{\rm AC}}
\newcommand\DS{\text{\rm{DS}}}

\newcommand\myStrut{\vphantom{\int^H}}

\title{
A Direct Mixed--Enriched Galerkin Method on Quadrilaterals
for Two-phase Darcy Flow
}

\author{Todd Arbogast \and Zhen Tao}

\institute{
  Todd Arbogast \at
  University of Texas at Austin;
  Department of Mathematics, C1200; Austin, TX 78712-1202 and
  Institute for Computational Engineering and Sciences, C0200;
  Austin, TX 78712-1229  \email{arbogast@ices.utexas.edu}
\and
 Zhen Tao\at
 University of Texas at Austin;
 Institute for Computational Engineering and Sciences, C0200;
 Austin, TX 78712--1229 \email{taozhen.cn@gmail.com}
}

\begin{document}

\maketitle

\begin{abstract}

We develop a locally conservative, finite element method for simulation of two-phase flow on quadrilateral meshes that minimize the number of degrees of freedom (DoFs) subject to accuracy requirements and the DoF continuity constraints.  We use a mixed finite element method (MFEM) for the flow problem and an enriched Galerkin method (EG) for the transport, stabilized with an entropy viscosity.  Standard elements for MFEM lose accuracy on quadrilaterals, so we use the newly developed AC elements which have our desired properties. Standard tensor product spaces used in EG have many excess DoFs, so we would like to use the minimal DoF serendipity elements.  However, the standard elements lose accuracy on quadrilaterals, so we use the newly developed direct serendipity elements.  We use the Hoteit-Firoozabadi formulation, which requires a capillary flux.  We compute this in a novel way that does not break down when one of the saturations degenerate to its residual value. Extension to three dimensions is described. Numerical tests show that accurate results are obtained.

\noindent

\noindent
\end{abstract}

\section{Introduction}\label{sec:Intro}

For many years now, finite element methods for Darcy flow and transport have been developed on rectangular and simplicial meshes. However, many methods lose accuracy when posed on distorted rectangular meshes.  In this work, we develop an accurate and efficient numerical method for two-phase flow in porous media on meshes of quadrilateral elements.  We also discuss extensions to meshes of cuboidal hexahedra for problems posed in three dimensions.

We are interested in non-rectangular meshes for at least three reasons: Firstly, space discretization on quadrilaterals uses half the number of elements compared to discretization on triangles with the same scale $h$, and on hexahedra, only $1/5$ or $1/6$ compared to tetrahedra.  Secondly, in most geoscience applications, discretization is determined by natural geologic layering. Lastly, problems with slightly deformable porous media require distorted meshes.

A classic procedure in the finite element method is to define appropriate finite elements on squares or cubes, which is relatively easy to do, and then to map these elements to quadrilaterals or hexahedra.  This works well under affine mappings.  However, a general quadrilateral or hexahedron is not affine equivalent to a reference square or cube.  Rather, the simplest map is bi- or tri-linear, which leads to inaccuracies in the approximation properties of the finite element spaces.  The fundamental idea we use to circumvent this problem, used also by other authors, is to define the finite element space directly on the physical element instead of mapping it from a reference element.  The direct finite element will contain polynomials, and so accurate approximation will result.

Two-phase flow in porous media is governed by a system of two partial differential equations (PDEs) combined with two equality constraints.  The PDEs can be formulated as a nearly elliptic or parabolic flow problem and a nearly hyperbolic transport problem.

Many locally conservative methods have been developed for elliptic flow problems.  We mention just a few of these, the mixed finite element method (MFE)~\cite{ERW_1984,Boffi_Brezzi_Fortin_2013}, the enhanced velocity method \cite{Wheeler_Wheeler_Yotov_2002,Thomas_Wheeler_2010}, the multipoint flux approximation methods (MPFA) and the multipoint flux mixed finite element methods (MFMFE) \cite{AEKWY_2007,WXY_2012,Wheeler_Yotov_2006_multipointFlux,Ingram_Wheeler_Yotov_2010}, the mimetic finite difference methods \cite{HMSS_2002}, and discontinuous Galerkin (DG) methods \cite{Wheeler_1978,BBO_1999,RWG_2001,Riviere_Wheeler_Girault_1999,Houston_Schwab_Suli_2002}.  We will approximate the flow problem using mixed finite element methods, which are locally conservative and give very accurate velocities.  However, as is well-known, classic mixed finite elements defined on a square or cube and mapped to a general convex quadrilateral or cuboidal hexahedron perform poorly; in fact, they fail to approximate the divergence in an optimal way (except the ABF and Devloo et al.\ spaces \cite{ABF_2002,Bergot_Durufle_2013,Siqueira_Devloo_Gomes_2013}, which have many excess DoFs). Recently, Arbogast and Correa \cite{Arbogast_Correa_2016} resolved the problem on quadrilaterals. They defined the AC spaces, two families of mixed finite elements that achieve optimal convergence properties on quadrilaterals, while maintaining efficiency by using the minimal number of Degrees of freedom (DoFs) possible. The current authors in \cite{Arbogast_Tao_2018_atSpaces} and Cockburn and Fu in \cite{Cockburn_Fu_2017_mDecompIII} defined H(div)-conforming mixed finite elements on cuboidal hexahedra that maintain accuracy and use the minimal number of DoFs on general hexahedra. (The current authors also defined some new mixed elements on quadrilaterals similar to the AC spaces in~\cite{Arbogast_Tao_2018x_serendipity}.)

Many methods have also been devised for approximation of the transport problem, including some of the methods mentioned above, ELLAM and characteristic methods \cite{Douglas_Russell_1982,ERW_1984,Arbogast_Wheeler_1995,WLELQ_2002,Arbogast_Huang_2006,WZET_2013}, the standard continuous Galerkin (CG) finite element method, and the enriched Galerkin (EG) method \cite{Sun_Liu_2009_eg,Lee_Lee_Wheeler_2016_eg}. The DG method utilizes discontinuous piecewise polynomial finite element spaces to approximate the solution and weakly enforces interelement continuity by penalty terms. However, the method uses a very high number of DoFs. The CG method uses many fewer DoFs, but it fails to provide locally conservative saturations, which can lead to non-physical results. We will use an EG method, which resolves this deficiency of the CG method.  Sun and Liu \cite{Sun_Liu_2009_eg} defined the EG method by enriching the approximation space of the CG method with elementwise constant functions and using it in the DG formulation. The EG method significantly reduces the global number of DoFs compared to DG, and, in fact, there is an efficient solver for elliptic and parabolic problems with EG approximations~\cite{Lee_Lee_Wheeler_2016_eg}.  The EG method has been applied to miscible displacement problems and two-phase flow in porous media in recent works~\cite{Lee_Wheeler_2018,Lee_Wheeler_2017_adaptive}.

We will base our finite element space for EG on serendipity finite elements \cite{Ciarlet_1978,Arnold_Awanou_2011}, since these use a fewer number of DoFs than full tensor product Lagrange finite elements.  However, it is well known that the accuracy of serendipity finite elements degrades when mapped from a reference element to a quadrilateral or hexahedron (at least for elements higher order than bilinear), while Lagrange finite elements maintain accuracy \cite{Arnold_Awanou_2011,Kaliakin_2001,Lee_Bathe_1993}.  Recently, the current authors introduced new, direct serendipity finite elements \cite{Arbogast_Tao_2018x_serendipity} that have the same number of degrees of freedom as the classic serendipity elements but maintain accuracy on general non-degenerate convex quadrilaterals.  We will use these elements.

There are many ways to separate the PDE system into flow and transport problems.  We use the approach promoted by Hoteit and Firoozabadi \cite{Hoteit_Firoozabadi_2008,Hoteit_Firoozabadi_2008_fractures}, solved using an IMPES operator splitting algorithm \cite{CHM_2006}.  The formulation requires construction of the divergence of the capillary flux, and we will provide a novel implementation that does not break down when the system degenerates (i.e., one of the saturations tends to its residual value). To avoid spurious oscillations due to sharp gradients in the solution of the transport problem, the stabilization technique of Guermond, Popov, and collaborators \cite{Guermond_Pasquetti_Popov_2011_entropyVisc,Bonito_Guermond_Popov_2014,ZGMP_2013} will be used.

To simplify the presentation, we concentrate on applying our new method to problems on quadrilateral meshes. In the next section we define the finite elements that we use on quadrilaterals.  The equations governing two-phase flow are given in Section~\ref{sec:twoPhase}.  We give the finite element method for flow in Section~\ref{sec:flow} and for transport in Section~\ref{sec:transport}. The IMPES coupling and extension to three dimensions is described briefly in Sections~\ref{sec:impes}--\ref{sec:3d}, respectively.  Numerical tests are given in Section~\ref{sec:numericalResults}, and we summarize and conclude our results in the final section.

\section{Direct finite elements on quadrilaterals}\label{sec:elem}

In this section, we review the finite element spaces that we will use.  Let $\Omega\subset\Re^2$ be a polygonal domain with boundary $\partial\Omega$, and let $\n$ denote the outward unit normal vector on the boundary.  We impose a conforming finite element mesh $\cT_h$ of quadrilaterals over the domain $\Omega$ of maximal diameter $h$. We will assume that the mesh $\cT_h$ is shape-regular \cite{Girault_Raviart_1986,Arbogast_Correa_2016}, which ensures that the mesh elements are not highly elongated nor degenerate nearly to triangles.  The outer unit normal to element $E\in\cT_h$ is $\n_{\partial E}$. The set of all interior edges or faces of $\cT_h$ is its skeleton, and it is denoted by $\cE_h$. Fix a unit normal vector $\n_\gamma$ for each edge $\gamma\in\cE_h$.

Let $\tilde\Po_k=\spn\{x^iy^{k-i}\,:\,i=0,1,\ldots,k\}$ be the vector space of homogeneous polynomials of exact degree~$k\ge0$, and let $\Po_k=\bigoplus_{n=0}^k\tilde\Po_n =\spn\{x^iy^{j}\,:\,i,j=0,1,\ldots,k;\ i+j\le k\}$ denote the space of polynomials of degree up to $k\ge0$.  Let $\Qu_{k}=\spn\{x^iy^{j}\,:\,i,j=0,1,\ldots,k\}$ be the space of tensor product polynomials of degree up to $k\ge0$. Finally, $\Po_k^2$ is the space of $2$-dimensional polynomial vectors for which each component is in $\Po_k$. We may specify that the domain of definition of the polynomials by writing, e.g., $\Po_k(E)$ for domain~$E$.

At times we will need a reference element, so fix it to be the square $\hat E=[-1,1]^2$.  For $E\in\cT_h$, the standard bilinear map $\F_E:\hat E\to E$ mapping vertices to vertices is bijective.  This gives rise to the map $\cF_E$, which maps a function $\hat f:\hat E\to\Re$ to a function $\cF_E(\hat f\,)=f:E\to\Re$ by the rule $f(\x) = \hat f(\hat\x)$, where $\x=\F_E(\hat\x)$.  We also have the Piola transform $\cP_E$ based on the bilinear map $\F_E$.  It maps a vector $\hat\v:\hat E\to\Re^2$ to a vector $\v:E\to\Re^2$, and it preserves the normal components $\hat\v\cdot\hat\n$ \cite{Boffi_Brezzi_Fortin_2013}.

\subsection{The AC mixed finite element spaces}\label{sec:AC}

Arbogast and Correa developed two families of mixed finite elements on quadrilateral meshes \cite{Arbogast_Correa_2016} for approximating $(\u,p)$ solving a second order elliptic equation in mixed form
\begin{equation}
\label{eq:strongEqn}
\u = -a\grad p,\quad\div\u = f\quad\text{in }\Omega,\quad\u\cdot\n = 0\quad\text{on }\partial\Omega,
\end{equation}
where $f\in L^2(\Omega)$ and the tensor $a$ is uniformly positive definite and bounded.  For $E\in\cT_h$ and index $s\ge0$, the AC elements are defined in terms of polynomials and a supplemental space of functions $\Su_s^\AC(E)$.  In terms of the vectors $\x\tilde\Po_s=\{(x,y)\,p(x,y):\:p\in\tilde\Po_s\}$, the family of full $H(\Div)$-approximation elements are
\begin{equation}
\label{eq:ac-div_E}
\V_\AC^s(E) = \Po_s^2\oplus\x\tilde\Po_s\oplus\Su_s^\AC(E)
\text{ and }
W_\AC^s(E) = \Po_s,
\end{equation}
while, the family of reduced $H(\Div)$-approximation elements are, for index $s\ge1$, 
\begin{equation}
\label{eq:ac-red_E}
\V_\AC^{s,\red}(E) = \Po_s^2\oplus\Su_s^\AC(E)
\text{ and }
W_\AC^{s,\red}(E) = \Po_{s-1}.
\end{equation}
The vector elements merge together to form $H(\Div)$ conforming spaces over $\Omega$ (i.e., the normal components of the vectors are continuous across each $\gamma\in\cE_h$), while the scalar elements remain discontinuous on $\Omega$. The former family approximates $\u$, $\div\u$, and $p$ to order $\cO(h^{s+1})$, while the latter family approximates $\u$ to order $\cO(h^{s+1})$ and $\div\u$ and $p$ or order $\cO(h^{s})$, as long as the mesh is shape regular.  These are the optimal convergence rates.  Moreover, the elements have the minimal number of DoFs possible for these rates under the restriction of $H(\Div)$ conformity.

Using the Piola transformation $\cP_E$, the space of supplemental vectors on $E$ is
\begin{equation}
\label{eq:Su}
\Su_s^\AC(E)=\begin{cases}
\spn\{\sigma_1^s,\sigma_2^s\},&s\ge 1,\\
\spn\{\sigma^0\},&s=0,
\end{cases}
\end{equation}
where $\sigma_i^s = \cP_E\hat\sigma_i^s$, $i=0,1,2$, and
\begin{equation}
\label{eq:hat_sigma}
\left\{\ \begin{alignedat}2
&\hat\sigma_1^s = \Curl(\hat x^{s-1}(1-\hat x^2)\hat y),&&\quad s\ge1,\\
&\hat\sigma_2^s = \Curl(\hat x\hat y^{s-1}(1-\hat y^2)),&&\quad s\ge1,\\
&\hat\sigma^0 = \Curl(\hat x\hat y),&&\quad s=0,
\end{alignedat}\right.
\end{equation}
are defined on the reference element $\hat E$.  Different from the implementation of classical mixed finite elements, the Piola mapping is only used to define the supplemental vectors. It is easy to implement these elements using the hybrid mixed finite element formulation \cite{Arnold_Brezzi_1985,Boffi_Brezzi_Fortin_2013}.


\subsection{Direct serendipity elements}\label{sec:DS}

\begin{figure}[htbp]\centering{\setlength\unitlength{3.2pt}
\begin{picture}(25,18)(-0.5,-0.5)
\thicklines
\put(5,16){\makebox(0,0){$r=2$}}
\put(0,0){\line(1,3){4}}
\put(0,0){\line(1,0){24}}
\put(24,0){\line(-1,4){4}}
\put(4,12){\line(4,1){16}}
\put(2,6){\circle*{1}}
\put(12,0){\circle*{1}}
\put(22,8){\circle*{1}}
\put(12,14){\circle*{1}}
\put(0,0){\circle*{1}}
\put(24,0){\circle*{1}}
\put(20,16){\circle*{1}}
\put(4,12){\circle*{1}}
\end{picture}\quad
\begin{picture}(25,18)(-0.5,-0.5)
\thicklines
\put(5,16){\makebox(0,0){$r=3$}}
\put(0,0){\line(1,3){4}}
\put(0,0){\line(1,0){24}}
\put(24,0){\line(-1,4){4}}
\put(4,12){\line(4,1){16}}
\put(1.333,4){\circle*{1}}
\put(2.667,8){\circle*{1}}
\put(8,0){\circle*{1}}
\put(16,0){\circle*{1}}
\put(22.667,5.333){\circle*{1}}
\put(21.333,10.667){\circle*{1}}
\put(9.333,13.333){\circle*{1}}
\put(14.667,14.667){\circle*{1}}
\put(0,0){\circle*{1}}
\put(24,0){\circle*{1}}
\put(20,16){\circle*{1}}
\put(4,12){\circle*{1}}
\end{picture}
}\centerline{\setlength\unitlength{3.2pt}
\begin{picture}(25,18)(-0.5,-0.5)
\thicklines
\put(5,16){\makebox(0,0){$r=4$}}
\put(0,0){\line(1,3){4}}
\put(0,0){\line(1,0){24}}
\put(24,0){\line(-1,4){4}}
\put(4,12){\line(4,1){16}}
\put(1,3){\circle*{1}}
\put(2,6){\circle*{1}}
\put(3,9){\circle*{1}}
\put(6,0){\circle*{1}}
\put(12,0){\circle*{1}}
\put(18,0){\circle*{1}}
\put(23,4){\circle*{1}}
\put(22,8){\circle*{1}}
\put(21,12){\circle*{1}}
\put(8,13){\circle*{1}}
\put(12,14){\circle*{1}}
\put(16,15){\circle*{1}}
\put(12,7){\circle*{1}}
\put(0,0){\circle*{1}}
\put(24,0){\circle*{1}}
\put(20,16){\circle*{1}}
\put(4,12){\circle*{1}}
\end{picture}\quad
\begin{picture}(25,18)(-0.5,-0.5)
\thicklines
\put(5,16){\makebox(0,0){$r=5$}}
\put(0,0){\line(1,3){4}}
\put(0,0){\line(1,0){24}}
\put(24,0){\line(-1,4){4}}
\put(4,12){\line(4,1){16}}
\put(0.8,2.4){\circle*{1}}
\put(1.6,4.8){\circle*{1}}
\put(2.4,7.2){\circle*{1}}
\put(3.2,9.6){\circle*{1}}
\put(4.8,0){\circle*{1}}
\put(9.6,0){\circle*{1}}
\put(14.4,0){\circle*{1}}
\put(19.2,0){\circle*{1}}
\put(20.8,12.8){\circle*{1}}
\put(21.6,9.6){\circle*{1}}
\put(22.4,6.4){\circle*{1}}
\put(23.2,3.2){\circle*{1}}
\put(7.2,12.8){\circle*{1}}
\put(10.4,13.6){\circle*{1}}
\put(13.6,14.4){\circle*{1}}
\put(16.8,15.2){\circle*{1}}
\put(12,9){\circle*{1}}
\put(9,6){\circle*{1}}
\put(15,6){\circle*{1}}
\put(0,0){\circle*{1}}
\put(24,0){\circle*{1}}
\put(20,16){\circle*{1}}
\put(4,12){\circle*{1}}
\end{picture}
}
\caption{ 
The nodal points for the DoFs of the direct serendipity finite element for small $r$.
\label{fig:nodalDoFs}}
\end{figure}
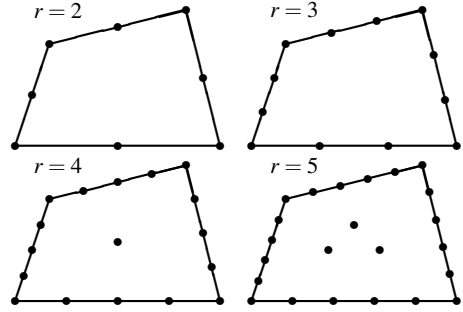

The current authors defined families of direct serendipity finite elements in a recent paper \cite{Arbogast_Tao_2018x_serendipity}. On $E\in\cT_h$, the element of index $r\ge1$ takes the form \begin{equation}\label{eq:serendipityForm}
\DS_{r}(E) = \Po_r\oplus\Su_r^\DS(E),
\end{equation}
and each family of elements is determined by the choice of the two supplemental functions spanning $\Su_r^\DS(E)$. The space can be merged into continuous (i.e., $H^1$ conforming) finite elements, and it has a minimal number of degrees of freedom (DoFs), which are depicted in Fig.~\ref{fig:nodalDoFs}. Its local dimension is
\begin{equation}
\dim\DS_{r}(E) = \dim\Po_r + 2 = \frac12(r+2)(r+1) + 2.
\end{equation}
A very general and explicit construction of these supplements is given in \cite{Arbogast_Tao_2018x_serendipity}.  They can be defined directly on $E$, or they can be defined on $\hat E$ and mapped to $E$.  In this work, we recall and use one of the simplest families of direct serendipity elements.

When $r=1$, it is known that the standard space of bilinear polynomials $\Qu_1$ defined on $\hat E$ and mapped to an element $E\in\cT_h$ maintains accuracy \cite{Arnold_Awanou_2011}.  These have the form \eqref{eq:serendipityForm}, with
\begin{equation}
\Su_1^\DS(E) = \spn\{\cF_E(\hat x\hat y)\}.
\end{equation}
Henceforth, we only describe the new direct serendipity finite elements for indices $r \geq 2$.

\begin{figure}[ht]\centering
\setlength\unitlength{3.2pt}
{\begin{picture}(34.5,32)(-5,-6)\small
%
\thicklines
\put(0,0){\line(1,3){4}}
\put(0,0){\line(1,0){24}}
\put(24,0){\line(-1,4){4}}
\put(4,12){\line(4,1){16}}
\put(2,6){\circle*{1}}\put(2,6){\vector(-3,1){4.74}}
\put(12,0){\circle*{1}}\put(12,0){\vector(0,-1){5}}
\put(22,8){\circle*{1}}\put(22,8){\vector(4,1){4.85}}
\put(12,14){\circle*{1}}\put(12,14){\vector(-1,4){1.21}}
\put(-4.3,7.5){\makebox(0,0){$\n_1$}}
\put(10.8,-5.9){\makebox(0,0){$\n_3$}}
\put(28.4,9){\makebox(0,0){$\n_2$}}
\put(10.8,20){\makebox(0,0){$\n_4$}}
\put(0,0){\circle*{1}}\put(-1.8,-1.6){\makebox(0,0){$\x_{v,13}$}}
\put(24,0){\circle*{1}}\put(25.6,-1.9){\makebox(0,0){$\x_{v,23}$}}
\put(20,16){\circle*{1}}\put(21.6,17.6){\makebox(0,0){$\x_{v,24}$}}
\put(4,12){\circle*{1}}\put(2.4,13.6){\makebox(0,0){$\x_{v,14}$}}
\put(12,7){\makebox(0,0){$E$}}
\put(4.6,8.8){\makebox(0,0){$e_1$}}
\put(8,-1.2){\makebox(0,0){$e_3$}}
\put(21.5,4){\makebox(0,0){$e_2$}}
\put(16,16.5){\makebox(0,0){$e_4$}}
\end{picture}}
\caption{A quadrilateral $E$, with edges $e_i$, outer unit normals $\n_i$, and vertices $\x_{v,13}$, $\x_{v,23}$, $\x_{v,24}$, and $\x_{v,14}$.
\label{fig:numbering}}
\end{figure}
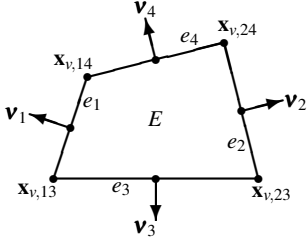

For the element $E\subset\cT_h$, let $\n_i$ denote the outer unit normal to edge $i$ (denoted $e_i$), $i=1,2,3,4$, and identify the vertices as $\x_{v,13}=e_1\cap e_3$, $\x_{v,14}=e_1\cap e_4$, $\x_{v,23}=e_2\cap e_3$, and $\x_{v,24}=e_2\cap e_4$ (see Fig.~\ref{fig:numbering}). We define the linear polynomial $\ell_i(\x)$ giving the distance of $\x\in\Re^2$ to edge~$e_i$ in the normal direction as
\begin{align}
\label{eq:lambda_i}
\ell_i(\x) &= - (\x-\x_i)\cdot\n_i, \quad i=1,2,3,4,
\end{align}
where $\x_i\in e_i$ is any point on the edge.  The functions are positive over the interior of $E$.

When $r\ge2$, we can define the shape functions associated to each vertex as
\begin{equation}\label{eq:shape-vertices}
\begin{alignedat}2
\phi_{v,13}(\x) &= \ell_2(\x)\ell_4(\x),&\quad\phi_{v,14}(\x) &= \ell_2(\x)\ell_3(\x),\\
\phi_{v,23}(\x) &= \ell_1(\x)\ell_4(\x),&\quad\phi_{v,24}(\x) &= \ell_1(\x)\ell_3(\x).
\end{alignedat}
\end{equation}
Interior shape functions appear when $r\ge4$. In this case, we take the shape functions
\begin{equation}\label{eq:shape-cell}
\phi_{E,j}\subset\ell_1\ell_2\ell_3\ell_4\Po_{r-4},\quad j=1,\ldots,\tfrac12(r-2)(r-3)
\end{equation}
so that they span the entire space $\ell_1\ell_2\ell_3\ell_4\Po_{r-4}$.  Since these are internal to $E$, the precise choice of basis is not so important.

The most interesting shape functions are those associated to the edges. Let
\begin{equation}\label{eq:simpleLambdaHV}
\ell_H = \ell_3 - \ell_4
\quad\text{and}\quad
\ell_V = \ell_1 - \ell_2,
\end{equation}
and define rational functions
\begin{align}
\label{eq:rational-v}
&R_V(\x) = \frac{\ell_1(\x) - \ell_2(\x)}{\xi_V^{-1}\ell_1(\x) + \eta_V^{-1}\ell_2(\x)},\\
\label{eq:rational-h}
&R_H(\x) = \frac{\ell_3(\x) - \ell_4(\x)}{\xi_H^{-1}\ell_3(\x) + \eta_H^{-1}\ell_4(\x)}
\end{align}
(note that the denominators do not vanish on $E$), where
$\n_H=(\n_3-\n_4)/|\n_3-\n_4|$, $\n_V=(\n_1-\n_2)/|\n_1-\n_2|$, and
\begin{align}
&\xi_V^{-1} = \sqrt{1-(\n_H\cdot\n_1)^2},\:
&\eta_V^{-1} = \sqrt{1-(\n_H\cdot\n_3)^2},\\
&\xi_H^{-1} = \sqrt{1-(\n_V\cdot\n_2)^2},\:
&\eta_H^{-1} = \sqrt{1-(\n_V\cdot\n_4)^2}.
\end{align}
Therefore,
\begin{alignat}3
\label{eq:RVe}
R_V(\x)|_{e_1} &= -\eta_V&&\quad\text{and}\quad&R_V(\x)|_{e_2} &= \xi_V,\\
\label{eq:RHe}
R_H(\x)|_{e_3} &= -\eta_H&&\quad\text{and}\quad&R_H(\x)|_{e_4} &= \xi_H.
\end{alignat}
There are $2(r-1)$ shape functions
associated to the edges $e_1$ and $e_2$, and they are
\begin{alignat}2
\label{eq:shape-h}
&\phi_{H,j}(\x) = \ell_3(\x)\ell_4(\x)\ell_H^{j}(\x),&&\quad j=0,1,\ldots,r-2,\\
\label{eq:shape-hpm}
&\phi_{H,r-1+j}(\x) =\ell_3(\x)\ell_4(\x)\rlap{$\ell_{V}(\x)\ell_H^{j}(\x),$}\\
\nonumber
&&&\quad j=0,1,\ldots,r-3,\\
&\label{eq:shape-hr}
\phi_{H,2r-3}(\x) =\ell_3(\x)\ell_4(\x)\rlap{$R_V(\x)\ell_H^{r-2}(\x).$}
\end{alignat}
In a similar way, shape functions associated with edges $e_3$ and $e_4$ are
\begin{alignat}2
\label{eq:shape-v}
&\phi_{V,j}(\x) = \ell_1(\x)\ell_2(\x)\ell_V^{j}(\x),&&\quad j=0,1,\ldots,r-2,\\
\label{eq:shape-vpm}
&\phi_{V,r-1+j}(\x) =\ell_1(\x)\ell_2(\x)\rlap{$\ell_{H}(\x)\ell_V^{j}(\x),$}\\
\nonumber
&&&\quad j=0,1,\ldots,r-3,\\
\label{eq:shape-vr}
&\phi_{V,2r-3}(\x) =\ell_1(\x)\ell_2(\x)\rlap{$R_H(\x)\ell_V^{r-2}(\x).$}
\end{alignat}
The edge shape functions are regular polynomials of degree~$r$ except the last two functions in each direction, which are rational functions.  However, all shape functions restrict to polynomials of degree~$r$ on the edges.

Finally,
\begin{equation}
\label{eq:supplementSpace}
\Su_r^{\DS}(E) = \spn\{\ell_3\ell_4\ell_H^{r-2}R_V,\ell_1\ell_2\ell_V^{r-2}R_H\},
\end{equation}
and
\begin{align}\label{eq:shapeSpace}
\DS_{r}(E) 
&= \spn\big\{ \phi_{v,13}, \phi_{v,14}, \phi_{v,23}, \phi_{v,24},\\
\nonumber
&\qquad\qquad\phi_{H,j}, \phi_{V,j}\ (j=0,1,\ldots,2r-3),\\
\nonumber
&\qquad\qquad\phi_{E,k}\ (k=1,\ldots,\tfrac12(r-2)(r-3))\big\}\\
\nonumber
&= \Po_r(E)\oplus\Su_r^\DS(E).
\end{align}
The unisolvence of the direct serendipity finite element space is proved generally in~\cite{Arbogast_Tao_2018x_serendipity}, or specifically for this choice of $R_V$, $R_H$, $\ell_V$ and $\ell_H$ in \cite{Arbogast_Tao_2017_serendipity}. A nodal basis can be constructed using local linear algebra as described in \cite{Arbogast_Tao_2017_serendipity,Arbogast_Tao_2018x_serendipity}.

\subsection{Enriched direct serendipity spaces}\label{sec:enrichedDS}

Let the discontinuous finite element spaces of order $r$
over $\cT_h$ be
\begin{align}
\label{eq:barXh-DS}
\bar X_r^{DS}(\cT_h) &= \left\{ \phi\in L^2(\Omega):\:\phi|_E \in \DS_r(E),\: E\in\cT_h\right\},\\
\label{eq:barXh-Q}
\bar X_r^{Q}(\cT_h) &= \left\{ \phi\in L^2(\Omega):\:\phi|_E=\cF_E(\hat\phi),\right.\\\nonumber
&\qquad\qquad\qquad\qquad\left.\hat\phi\in\Qu_r(\hat E),\: E\in\cT_h\right\},\\
\label{eq:barXh-P}
\bar X_r^{P}(\cT_h) &= \left\{ \phi\in L^2(\Omega):\:\phi|_E \in \Po_r(E),\: E\in\cT_h\right\}.
\end{align}
We remark that $\bar X_0^{Q} = \bar X_0^{P}$ are the spaces of piecewise constants on each element. The continuous finite element spaces of order $r$ over $\cT_h$ are
\begin{align}
\tilde X_r^{DS}(\cT_h) &= \bar X_r^{DS}(\cT_h) \cap \cC^0(\Omega),\\ 
\tilde X_r^{Q}(\cT_h) &= \bar X_r^{Q}(\cT_h) \cap \cC^0(\Omega),
\end{align}
where $\cC^0(\Omega)$ is the set of continuous functions over $\Omega$.
Finally, we define the enriched finite element spaces of order $r$
over $\cT_h$ as
\begin{align}
\label{eq:Xh-DS}
 X_r^{DS}(\cT_h) &= \tilde X_r^{DS}(\cT_h) + \bar X_0^{Q}(\cT_h), \\
\label{eq:Xh-Q}
 X_r^{Q}(\cT_h) &= \tilde X_r^{Q}(\cT_h) + \bar X_0^{Q}(\cT_h).
\end{align}

When used in a standard discontinuous Galerkin weak formulation, we can obtain locally mass conservative, discontinous or enriched Galerkin approximations.  We denote these methods by DG-$\Qu_r$, DG-$\Po_r$, EG-$\Qu_r$ and EG-$\DS_r$, using, respectively, the elements $\bar X_r^{Q}(\cT_h)$, $\bar X_r^{P}(\cT_h)$, $X_r^{Q}(\cT_h)$, and $X_r^{DS}(\cT_h)$.  In Table~\ref{tab:eg:compare-2d}, we show the number of degrees of freedom of these four methods for an $n\times n$ two dimensional quadrilateral mesh. All the methods achieve the same convergence rate.  The EG-$\DS_r$ methods ultilize the fewest number of degrees of freedom.

\begin{table}[ht]
\caption{A comparison of the global number of DoFs of
DG-$\Qu_r$, DG-$\Po_r$, EG-$\Qu_r$, and EG-$\DS_r$ on an $n\times n$ mesh.}
\label{tab:eg:compare-2d}
\centering
\begin{tabular}{c|cccc}
\hline
$\myStrut r$ & DG-$\Qu_r$ & DG-$\Po_r$ & EG-$\Qu_r$ & EG-$\DS_r$\\
\hline
$\myStrut 0$ & $n^2$ & $n^2$ & --- & ---\\
$1$ & $4n^2$ & $3n^2$ & $2n^2+2n+1$ & $2n^2+2n+1$ \\
$2$ & $9n^2$ & $6n^2$ & $5n^2+4n+1$ & $4n^2+4n+1$ \\
$3$ & $16n^2$ & $10n^2$ & $10n^2+6n+1$ & $6n^2+6n+1$ \\
$4$ & $25n^2$ & $15n^2$ & $17n^2+8n+1$ & $9n^2+8n+1$ \\
\hline 
\end{tabular}
\end{table}

\section{Two-phase flow formulation}\label{sec:twoPhase}

Let subscript $w$ denote the wetting phase and subscript $n$ the non-wetting phase.  The pressures, Darcy velocities, and saturations of each phase are denoted $p_\alpha$, $\u_\alpha$, and $S_\alpha$, respectively, $\alpha = w,n$. The Darcy velocity of each phase satisfies
\begin{equation}\label{eq:2p:darcy-law}
\u_\alpha = - \frac{k_{r\alpha}(S_{w})}{\mu_\alpha}\mathbf{K} (\grad p_\alpha - \rho_\alpha g\grad z),
\quad \alpha = w,n,
\end{equation}
where $k_{r\alpha}$ is the relative permeability of the $\alpha$ phase, which depends on the water saturation $S_w$, $\mathbf{K}$ is the absolute rock permeability, $\mu_\alpha$ and $\rho_\alpha$ are the viscosity and density of the $\alpha$ phase (here assumed constant), $g$ is the gravitational constant, and $\grad z$ defines the direction of gravity. Volume balance requires the algebraic constraint
\begin{equation}\label{eq:2p:sat-conservation}
S_w + S_n = 1,
\end{equation}
and conservation of mass requires
\begin{equation}\label{eq:2p:mass-conservation} 
\phi\frac{\partial S_\alpha}{\partial t} + \div\u_\alpha = q_\alpha, \quad \alpha = w,n,
\end{equation}
where $\phi$ is the porosity and $q_\alpha$ are the source or sink (i.e., well) terms.
One also needs the capillary pressure relation
\begin{equation}\label{eq:2p:pc0}
p_c(S_w) = p_n - p_w.
\end{equation}

One typical capillary pressure function is given by
\begin{align}\label{eq:2p:pc}
p_c(S_w) = -B_c \log(S_e),
\end{align}
where $B_c$ is a positive parameter inversely proportional to $\sqrt{k}$ and the normalized saturation $S_e$ is given by
\begin{align}\label{eq:2p:se}
S_e = \frac{S_w - S_{rw}}{1-S_{rw}-S_{rn}},
\end{align}
where $S_{r\alpha}$ are the residual saturations for the wetting and non-wetting phases.
An simple example of the relative permeability of each phase is given by
\begin{align}\label{eq:2p:krs}
k_{rw} = S_e^\beta\quad\text{and}\quad k_{rn}=(1-S_e)^\beta,
\end{align}
where $\beta$ is a positive parameter. Other commonly used rock models are those of Brooks-Corey and van Genuchten \cite{Bear_1972,CHM_2006}.

For a classical formulation, introduce the phase mobilities
\begin{equation}
\lambda_\alpha = \frac{k_{r\alpha}}{\mu_\alpha}, \qquad \alpha = w,n,
\end{equation}
the total mobility
\begin{equation}
\lambda_t = \lambda_w + \lambda_n,
\end{equation}
and the fractional flow functions
\begin{equation}
f_\alpha = \frac{\lambda_\alpha}{\lambda_t}, \qquad \alpha = w,n.
\end{equation}
Sum the mass conservation equations \eqref{eq:2p:mass-conservation} and combine with the algebraic 
constraint \eqref{eq:2p:sat-conservation} to obtain the \emph{flow equation}
\begin{align}\label{eq:2p:ut}
\div\u_t = q_t,
\end{align}
where $\u_t=\u_w + \u_n$ is the total flow velocity, and
\begin{align}\label{eq:2p:ut2}
\u_t =  -\lambda_t\mathbf{K}\grad p_w - \lambda_n\mathbf{K}\grad p_c
 + (\rho\lambda)_t g\mathbf{K}\grad z,
\end{align}
where $q_t=q_w+q_n$, $(\rho\lambda)_t = \rho_w\lambda_w+\rho_n\lambda_n$.
Therefore, the phase velocity $\u_w$ and $\u_n$ are related to the total velocity by
\begin{align}\label{eq:2p:water-vel}
\u_w &= f_w \u + \mathbf{K}\lambda_n f_w \grad p_c +\mathbf{K} \lambda_n f_w(\rho_w - \rho_n)g\grad z,\\
\u_n &= f_n \u - \mathbf{K}\lambda_w f_n \grad p_c +\mathbf{K} \lambda_w f_n(\rho_n - \rho_w)g\grad z.
\end{align}
The \emph{saturation equation} for the water phase becomes
\begin{align}\label{eq:2p:saturation}
&\phi \frac{\partial S_w}{\partial t} + 
\div \Big\{ \mathbf{K}f_w(S_w)\lambda_n(S_w)\Big(\frac{dp_c}{dS_w}\grad S_w \\\nonumber
&\qquad\qquad\qquad + (\rho_w - \rho_n)g\grad z\Big)
+f_w(S_w)\u \Big\} =  q_w.
\end{align}
We have three unknowns $\{p_w, \u_t, S_w\}$ in the three equations \eqref{eq:2p:ut}, \eqref{eq:2p:ut2}, and \eqref{eq:2p:saturation}, which together comprise a well defined system (see, e.g., \cite{Arbogast_1992a}).

The classical formulation is challenged when the capillary pressure is degenerate or discontinuous.  Suppose the water saturation is near the wetting phase residual saturation $S_{rw}$. Then the non-wetting phase mobility is about ${1}/{\mu_n}$, and at the same time the capillary pressure and its derivative tends to infinity, i.e., $p_c(S_{rw})\to\infty$ and $p_c'(S_{rw})\rightarrow\infty$.  Thus the term $\lambda_n\mathbf{K}\grad p_c$ in~\eqref{eq:2p:ut2} also tends to infinity.  If we choose the non-wetting phase pressure $p_n$ as our primary variable instead of $p_w$, we can avoid this singularity around $S_{rw}$. However, if the derivative of the capillary pressure also tends to infinity at the non-wetting phase residual saturation, e.g., in a van Genuchten model, we will meet the same problem even with the non-wetting phase pressure as the primary variable.  Moreover, when there are different rock types present, i.e., capillarity is heterogeneous as in \eqref{eq:2p:pc}, the spatial gradient of the capillary pressure $\grad p_c$ also tends to infinity.

\subsection{The Hoteit-Firoozabadi formulation}

In~\cite{Hoteit_Firoozabadi_2008,Hoteit_Firoozabadi_2008_fractures}, Hoteit and Firoozabadi presented a new formulation that avoids the drawbacks of the classical one. They use the classic flow potential
\begin{align}\label{eq:2p:phi}
\Phi_\alpha = p_\alpha+\rho_\alpha gz, \quad \alpha=w,n,
\end{align}
which leads to the capillary potential
\begin{align}
\Phi_c = \Phi_n-\Phi_w = p_c+(\rho_n-\rho_w)gz.
\end{align}
The total velocity $\u_t$ is then written in terms of two velocity variables
$\u_a$ and $\u_c$, as follows:
\begin{align}
\u_t & = \u_a+\u_c,\\
\label{eq:2p:ua}
\u_a & = -\lambda_t\mathbf{K}\grad\Phi_w,\\
\label{eq:2p:uc}
\u_c &= -\lambda_n\mathbf{K}\grad\Phi_c.
\end{align}

The velocity variable $\u_a$ has the same pressure driving force as the wetting phase velocity $\u_w$,
but it has a smoother mobility $\lambda_t$ rather than $\lambda_w$. Since the total mobility $\lambda_t$ is strictly positive, the value of $\lambda_t^{-1}$ is bounded. The wetting phase velocity is then simply
\begin{align}
\u_w = \frac{\lambda_w}{\lambda_t}(-\lambda_t\mathbf{K}\grad \Phi_w) = f_w\u_a.
\end{align}
We rewrite the flow equation with the new splitting and the saturation equation in terms of $\u_a$ to obtain
\begin{align}
\label{eq:2p:new-flow}
&\div\u_a = q_t-\div\u_c,\\
\label{eq:2p:new-sat}
&\phi\frac{\partial S_w}{\partial t} = q_w -\div(f_w\u_a).
\end{align}

\section{Approximation of flow and construction of the capillary flux}\label{sec:flow}

For ease of exposition, we collect the equations from the previous section that govern flow, and we add appropriate boundary conditions.  At a fixed time~$t$,
\begin{alignat}2
\label{eq:flow_darcy}
\u_a &= -\lambda_t\mathbf{K}\grad\Phi_w &&\quad\text{in }\Omega,\\
\label{eq:flow_conservation}
\div\u_a &= q_t-\div\u_c &&\quad\text{in }\Omega,\\
\label{eq:flow_bc_pressure}
\Phi_w &= \hat\Phi_B &&\quad\text{on }\Gamma_D,\\
\label{eq:flow_bc_flux}
\u_a\cdot\n &=  u_B-\u_c\cdot\n &&\quad\text{on }\Gamma_N,
\end{alignat}
where $\partial\Omega = \Gamma_N\cup\Gamma_D$ has been decomposed into disjoint Neumann and Dirichlet parts of the boundary and $\hat\Phi_B$ and $u_B$ are given.  We assume that $\gamma\in\cE_h$ lies within either $\Gamma_D$ or $\Gamma_N$, and denote $\cE_{h,D}=\cE_h\cap\Gamma_D$ and $\cE_{h,N}=\cE_h\cap\Gamma_N$. At each time step, we will solve the system for the unknowns $\u_a$ and $\Phi_w$, assuming that $\u_{c}$ and the wetting phase saturation $S_{w}$ are given.

\subsection{A mixed method for the flow equation}

For $E\in\cT_h$, let $\V_h(E)\times W_h(E)$ denote the mixed finite elements used, which in our case are AC elements of index $s\geq0$, as described in Section~\ref{sec:AC} (although one could use any of the direct mixed finite elements defined in \cite{Arbogast_Tao_2018x_serendipity}).  If one restricted to rectangular meshes, one could use the classic Raviart-Thomas (RT) elements \cite{Raviart_Thomas_1977}.  We solve \eqref{eq:flow_darcy}--\eqref{eq:flow_bc_flux} using the hybrid form of the mixed method \cite{Arnold_Brezzi_1985,Boffi_Brezzi_Fortin_2013}.  To this end, we define the Lagrange multiplier space $M_h$ to be the set of piecewise polynomials of degree up to $s$ defined on each of the skeleton $\cE_h$. Then the global mixed spaces are 
\begin{align*}
\V_h &=\{\v_h:\v_h|_E\in\V_h(E)\ \forall E\in\cT_h\},\\
W_h &=\{w_h:w_h|_E\in W_h(E)\ \forall E\in\cT_h\},\\
M_h &= \{\mu_h:\mu_h|_\gamma\in \Po_s(\gamma)\ \forall \gamma\in\cE_h, \gamma\not\in\Gamma_D\}.
\end{align*}
Note that we do not enforce $H(\Div)$ conformity on $\V_h$.

The hybrid mixed finite element formulation is: Find $\u_{a,h}\in \V_h$, $\Phi_{w,h}\in W_h$,
and $\hat\Phi_{w,h}\in M_h$ such that
\begin{align}
\label{eq:2p:flow1}
&\int_\Omega (\lambda_t\mathbf{K})^{-1}\u_{a,h}\v_h
- \sum_{E\in\cT_h}\int_E \Phi_{w,h}\div\v_h\\ \nonumber
&\quad + \sum_{E\in\cT_h}\int_{\partial E\backslash\Gamma_D}\!\hat\Phi_{w,h} \v_h\cdot\n_{\partial E}
= -\int_{\Gamma_D}\!\hat\Phi_{B} \v_h\cdot\n\ \ \forall\v_h\in\V_h,\\
\label{eq:2p:qc}
&\sum_{E\in\cT_h}\int_E \div\u_{a,h}\,w_h = \int_\Omega (q_t-\div\u_{c,h})\,w_h\quad\forall w_h\in W_h,\\
\label{eq:2p:flow3}
&\sum_{E\in\cT_h}\int_{\partial E\backslash\Gamma_D}\!\u_{a,h}\cdot\n_{\partial E}\,\mu_h
 = \int_{\Gamma_N}\!(u_B-\u_c\cdot\n)\,\mu_h\\[-8pt] \nonumber
&\qquad\qquad\qquad\qquad\qquad\qquad\qquad\qquad\forall \mu_h\in M_h.
\end{align}
Since $\lambda_t^{-1}$ is always bounded, the first term in \eqref{eq:2p:flow1} is well-defined. The approximation spaces $\V_h$ and $W_h$ are defined elementwise and the continuity is only enforced through the Lagrange multiplier $\hat\Phi_w$ on the skeleton of the mesh $\cT_h$. In terms of the DoFs of the solution, we obtain the following algebraic problem 
\begin{equation}\label{eq:2p:abl}
\begin{pmatrix}
A_\lambda  & \; -B  & \; L\\
B^T & \; 0 & \; 0\\
L^T & \; 0 & \; 0
\end{pmatrix}
\begin{pmatrix}
\u_{a,h} \\
\Phi_{w,h} \\
\hat\Phi_{w,h}
\end{pmatrix}=
\begin{pmatrix}
\Phi_B\\
Q \\
U_B
\end{pmatrix},
\end{equation}
where $A_\lambda$ is symmetric positive definite and block diagonal, with block size equal to $\dim\V_h(E)$. Using the Schur complement technique, we can express
\begin{align}\label{eq:2p:hf17}
\u_{a,h} = A_\lambda^{-1}\big(\Phi_B + B\Phi_{w,h} - L\hat\Phi_{w,h}\big),
\end{align}
and therefore
\begin{align}
\label{eq:2p:hf20}
L^{T}\u_{a,h} &= L^{T}A_\lambda^{-1}\big(\Phi_B + B\Phi_{w,h} - L\hat\Phi_{w,h}\big) = Q,\\
\label{eq:2p:hf25}
B^{T}\u_{a,h} &= B^TA_\lambda^{-1}\big(\Phi_B + B\Phi_{w,h} - L\hat\Phi_{w,h}\big) = U_B.
\end{align}
If the lowest order RT spaces are used, \eqref{eq:2p:hf17}--\eqref{eq:2p:hf25} is equivalent to equations used
in the work of Hoteit and Firoozabadi \cite[equations (17), (20), (25)]{Hoteit_Firoozabadi_2008}.

\subsection{Construction of the capillary flux}

The capillary flux $\u_c$ is defined by \eqref{eq:2p:uc}, i.e., $\u_{c} = -\lambda_n\mathbf{K}\grad \Phi_c$. We have assumed the saturation $S_w$ is given, so $\Phi_c(S_w)$ is known inside each element.  However, $S_w$ is discontinuous on $\gamma\in\cE_h$, so $\lambda_n$ is not known on the skeleton.  Moreover, it vanishes at $S_w=1-S_{rn}$, so we cannot invert it.  A weak form of the equation is to solve it locally on each $E\in\cT_h$ as follows: Find
$\u_{c,h}\in\V_h$ and $\hat\Phi_{c,h}\in M_h$ such that
\begin{align}\label{eq:cap-flux-weak}
&\int_E\mathbf{K}^{-1} \u_{c,h}\v_h + \int_{\partial E}\hat\Phi_{c,h}\,\hat\lambda_n\v_h\cdot\n_\gamma\\\nonumber
&\quad = \int_E\Phi_{c}(S_w)\div\left(\lambda_n(S_w)\v_h\right)\quad\forall\v_h\in\V_h(E),
\end{align}
and impose continuity of the capillary flux
\begin{equation}\label{eq:cap-flux-weak-continuity}
\sum_{E\in\cT_h}\int_{\partial E} \u_{c,h}\cdot\n_\gamma \mu_h=0\quad\forall\mu_h\in M_h,
\end{equation}
where $\hat\lambda_n$ is the interface value of $\lambda_n(S_w)$.  Note that in this model, the capillary flux is set to zero on $\partial\Omega$.

Following \cite{Hoteit_Firoozabadi_2008}, the matrix form of the equations can be written as
\begin{align}
\begin{pmatrix}
A  & \; L_\lambda\\
L^T & \; 0
\end{pmatrix}
\begin{pmatrix}
\u_{c,h} \\
\hat\Phi_{c,h}
\end{pmatrix}=
\begin{pmatrix}
B_\lambda\Phi_c\\
0
\end{pmatrix},
\end{align}
and we can eliminate the variable $\u_{c,h}$ through a Schur complement, i.e.,
\begin{align}
\label{eq:2p:uc-rec}
\u_{c,h} &= A^{-1}\big(B_\lambda\Phi_c - L_\lambda\hat\Phi_{c,h}\big),\\
\label{eq:2p:hf31}
L^T\u_{c,h} &= L^TA^{-1}\big(B_\lambda\Phi_c - L_\lambda\hat\Phi_{c,h}\big) = 0.
\end{align}
Knowing $\Phi_c$, we solve the latter equation for $\hat\Phi_{c,h}$ and construct $\u_{c,h}$ from the former equation. Equation~\eqref{eq:2p:hf31} is equivalent to \cite[equation (31)]{Hoteit_Firoozabadi_2008}. The matrix $L^TA^{-1}L_\lambda$ is singular if $\lambda_n=0$, and in that case, the linear system~\eqref{eq:2p:hf31} cannot be solved. In~\cite{Hoteit_Firoozabadi_2008}, an upstream strategy is applied to resolve the singularity. However, if the whole region is filled with the wetting phase ($S_e=1$), we can not borrow a nonzero $\lambda_n$ from any neighboring elements.

In order to solve the degenerate equation, we provide an alternate approach.  Let $\hat\zeta_h = \hat\lambda_n\hat\Phi_{c,h}$, and solve in place of \eqref{eq:cap-flux-weak}: Find $\u_{c,h}\in \V_h$ and $\hat\zeta_h\in M_h$ such that
\begin{align}
\label{eq:2p:uc-re1}
&\int_E\mathbf{K}^{-1} \u_{c,h}\v_h + \int_{\partial E}\hat\zeta_h\v_h\cdot\n_\gamma\\\nonumber
&\quad = \int_E\Phi_{c}(S_w)\div\left(\lambda_n(S_w)\v_h\right) \quad\forall\v_h\in \V_h(E).
\end{align}
Now the matrix form of this with \eqref{eq:cap-flux-weak-continuity} is
\begin{align}
\begin{pmatrix}
A  & \; L\\
L^T & \; 0
\end{pmatrix}
\begin{pmatrix}
\u_{c,h} \\
\hat\zeta_h
\end{pmatrix}=
\begin{pmatrix}
B_\lambda\Phi_c\\
0
\end{pmatrix},
\end{align}
and so
\begin{align}
\label{eq:2p:uc-rec2}
\u_{c,h} &= A^{-1}\big(B_\lambda\Phi_c - L\hat\zeta_h\big),\\
\label{eq:2p:hf312}
L^T\u_{c,h} &= L^TA^{-1}\big(B_\lambda\Phi_c - L\hat\zeta_h\big) = 0.
\end{align}
The Schur complement $L^TA^{-1}L$ is now symmetric and does not depend on the value of saturation. The solution can be shifted by a constant and the differential equation still holds, so we need to set a constraint on the average or assign a Dirichlet boundary condition on a small portion of the boundary. Then the Schur complement system is non-singular and can always be solved. Once $\hat\zeta_h$ is found, we can construct $\u_{c,h}$ by~\eqref{eq:2p:uc-rec2} and finally we obtain the capillary source $B^T\u_{c,h}$ needed in \eqref{eq:2p:qc}.

We remark that in~\cite{Lee_Wheeler_2018} and~\cite{AJPW_2013}, the authors enforce the continuity of the capillary pressure $p_c$ (or $\Phi_c$ with the gravity igonored) by penalizing the interface jump in $p_c$. In this work and in \cite{Hoteit_Firoozabadi_2008,Hoteit_Firoozabadi_2008_fractures}, weak continuity of the capillary flux $\u_c$ is enforced across the face, and no penalty parameter is required.

\section{Approximation of transport with entropy stabilization}\label{sec:transport}

The saturation equation \eqref{eq:2p:new-sat} will be solved assuming that $\u_a$ is known.  In that case, it is of hyperbolic type for the saturation. For simplicity, we impose Neumann conditions on the boundary $\partial\Omega$. To this end, decompose it into two disjoint parts, $\partial\Omega=\Gamma_{in}\cup\Gamma_{out}$, where the inflow boundary is $\Gamma_{in}=\left\{\x\in\partial\Omega:\:\u_a\cdot\n<0 \right\}$, and the outflow boundary is $\Gamma_{out}=\left\{\x\in\partial\Omega:\:\u_a\cdot\n\geq 0 \right\}$.  We impose the inflow boundary condition $S_w=S_B$ on $\Gamma_{in}$, where $S_B$ is given.

Because \eqref{eq:2p:new-sat} is hyperbolic, numerical schemes need to be stabilized.  We choose to stabilize our hyperbolic equation by adding an entropy viscosity $\mu_h$. The stabilized transport problem, posed for $t>0$, becomes
\begin{alignat}2
\label{eq:2p:new-sat2}
&\phi\frac{\partial S_w}{\partial t} -\div(\mu_h\grad S_w) =q_w -\div(f_w\u_a)&&\quad\text{in }\Omega,\\
\label{eq:stabTranspBCin}
&(S_w\u_a-\mu_h\grad S_w)\cdot\n = S_B\,\u_a\cdot\n &&\quad\text{on }\Gamma_{in},\\
&-\mu_h\grad S_w\cdot\n = 0 &&\quad\text{on }\Gamma_{out},\\
\label{eq:2p:new-ic}
&S_w(\x,0) = S_0(\x)&&\quad\text{for }\x\in\rlap{$\Omega$.}
\end{alignat}
The inflow condition \eqref{eq:stabTranspBCin} corresponds to $S_w = S_B$ for the unstabilized equation \eqref{eq:2p:new-sat} (i.e., $\mu_h=0$).

\subsection{Enriched Galerkin for the saturation equation}

Let subscript $m=0,1,\ldots$ denote the time step index.  Within a transport time step $(t^m,t^{m+1}]$, we assume that we know the velocity $\u_a^m$ and the wetting phase saturation $S_w^m$ at the earlier time~$t^m$.

We assume that the mesh skeleton $\cE_h$ is aligned with the interfaces between $\Gamma_{in}$ and $\Gamma_{out}$.  Let the set of boundary faces in $\cE_h$ be denoted by $\Gamma_{h,in}$ and $\Gamma_{h,out}$.  Suppose two elements $E_i$ and $E_j$ in $\cT_h$ are neighbors, so $\gamma=\partial E_i\cap\partial E_j\in\cE_h$ is nonempty, and that the unit normal $\n_\gamma$ points from $E_i$ to $E_j$. For a function $\phi$ which is piecewise continuous on each $E\in\cT_h$, we define the average and jump on $\gamma$ as
\begin{align}
\avg{\phi}_\gamma &= \dfrac12 \left( (\phi|_{E_i})\big|_\gamma+(\phi|_{E_j})\big|_\gamma \right),\\
\jmp{\phi}_\gamma &=  (\phi|_{E_j})\big|_\gamma - (\phi|_{E_i})\big|_\gamma.
\end{align}
Let the upwind value of the interface saturation be
\begin{align}
S_w^{*,m}|_\gamma = 
\begin{cases}
S_w^m|_{E_i}\quad \text{if }\u_a^m\cdot\n_\gamma\geq 0,\\
S_w^m|_{E_j}\quad \text{if }\u_a^m\cdot\n_\gamma < 0.
\end{cases}
\end{align}

Define the interior penalty term $J_{\sigma}$ mapping to $\Re$ as
\begin{align}
\label{eq:eg:j1}
J_{\sigma}(s,w) &= \sum_{\gamma\in\cE_h}
\dfrac{r^2\sigma_\gamma}{|\gamma|^{1/(d-1)}}\int_\gamma\,\jmp{s}_\gamma\jmp{w}_\gamma,
\end{align}
where $\sigma_\gamma>0$ is the penalty parameter for face $\gamma\in\cE_h$, $r$ is the degree of the polynomials used in the finite element space, $d=2$ is the dimension of the domain, and $|\gamma|$ is the measure of $\gamma$. Let
\begin{align}\label{eq:eg:a2-bilinear}
a_{\text{DG}}^m(s,w) = & \sum_{E\in\cT_h}\int_E \mu_h^m\grad s\cdot\grad w + J_\sigma(s,w)\\\nonumber
& -\sum_{\gamma\in\cE_h}\int_\gamma\avg{\mu_h^m\grad s\cdot\n_\gamma}_\gamma\jmp{w}_\gamma \\\nonumber
& - s_{\text{form}}\sum_{\gamma\in\cE_h}\int_\gamma\avg{\mu_h^m\grad w\cdot\n_\gamma}_\gamma\jmp{s}_\gamma,
\end{align}
where $\mu_h^m$ will be defined in the next subsection and $s_{\text{form}}\in\{-1,0,1\}$ determines whether one chooses to use the non-symmetric (NIPG, $s_{\text{form}}=-1$), incomplete (IIPG, $s_{\text{form}}=0$), or symmetric (SIPG, $s_{\text{form}}=1$) interior penalty Galerkin formulation. Let the linear functional $L_{\text{DG}}^m(w)$ for the time step $m$ be
\begin{align}
L_{\text{DG}}^m(w) &= \int_\Omega q_w^m w+ \sum_{E\in\cT_h}\int_E f_w(S_w^m)\,\u_a^m\cdot\grad w\\
&\quad -\sum_{\gamma\in\cE_h}\int_\gamma f_w(S_w^{*,m})\u_a^m\cdot\n_\gamma\jmp{w}_\gamma\nonumber\\
&\quad-\sum_{\gamma\in\Gamma_{h,out}^m}\int_\gamma f_w(S_w^m)\u_a^m\cdot\n_\gamma w \nonumber\\
&\quad -\sum_{\gamma\in\Gamma_{h,in}^m}\int_\gamma f_w(S_B^m)\u_a^m\cdot\n_\gamma w.\nonumber
\end{align}
Finally, after applying backward Euler for the time discretization, we have the fully
discrete Galerkin formulation: Given $S_{w,h}^m$, find $S_{w,h}^{m+1}\in X_h(\cT_h)$ such that
\begin{alignat}2
\label{eq:2p:trans1}
&\int_\Omega\dfrac{\phi S_{w,h}^{m+1}-\phi S_{w,h}^m}{\delta t^{m+1}}\,w_h
+ \rlap{$a_{\text{DG}}^m(S_{w,h}^{m+1}, w_h)$}\\\nonumber
&\qquad  = L_{\text{DG}}^m(w_h)&&\quad\forall w_h\in X_h(\cT_h),\\
\label{eq:2p:trans1-ic}
& (\phi S_{w,h}^0, w_h) = (\phi S_0, w_h)&&\quad\forall w_h\in X_h(\cT_h),
\end{alignat}
where $\delta t^{m+1}=t^{m+1}-t^m$ and $X_h(\cT_h)$ can be any of the discontinuous Galerkin finite element spaces defined in \eqref{eq:barXh-DS}--\eqref{eq:barXh-P} or any of the enriched spaces defined in \eqref{eq:Xh-DS}--\eqref{eq:Xh-Q}.

\subsection{Entropy viscosity}

We review the entropy viscosity method of Guermond, Pasquetti, and Popov \cite{Guermond_Pasquetti_Popov_2011_entropyVisc} for the conservation law
\begin{equation}
\partial_t c + \div (F(c)) = 0.\label{eq:eg:cons-law}
\end{equation}
An entropy pair is a set of functions $\En(c)$ and $\mathsf{F}(c)$ such that $\En(c)$ is convex and $\mathsf{F}'(c) = \En'(c)F'(c)$. The function $\En(c)$ is called an entropy, and $\mathsf{F}(c)$ is called the associated entropy flux. A commonly used choice is $\En(c) = \frac{1}{2}c^2$.  Then the entropy solution of \eqref{eq:eg:cons-law} satisfies
\begin{equation}
\partial_t \En(c) + \div (\mathsf{F}(c)) \leq 0,
\end{equation}
with equality holding except at the shocks.

Given a numerical approximation $c_h(\cdot,t)$ at time $t$, we define the entropy residual
\begin{align}
\label{eq:eg:ent-res}
R_h(\cdot,t) = \partial_t \En(c_h(\cdot,t)) + \div\mathsf{F}(c_h(\cdot,t)).
\end{align}
Following~\cite{Guermond_Pasquetti_Popov_2011_entropyVisc,ZGMP_2013,Bonito_Guermond_Popov_2014,Lee_Wheeler_2017_adaptive}, we use this residual to define a viscosity $\mu_\En$ on each $E\in\cT_h$ as
\begin{align}
\label{eq:eg:nu-e}
\mu_\En = \lambda_\En\,h_E^2\,\dfrac{\|R_h\|_{\infty,E}}{\|\En(c_h)-\bar{\En}(c_h)\|_{\infty,\Omega}},
\end{align}
where $\lambda_\En$ is a tunable parameter, $h_E=\text{diam}(E)$, $\|\cdot\|_{\infty}$ is the maximum norm, the denominator is a normalization constant, and $\bar{\En}(c_h) = \frac{1}{|\Omega|}\int_\Omega\En(c_h)$. However, in our applications, the (true) saturation is always between 0 and 1, and we intend to take the entropy $\En(c)=\frac12 c^2$, so the normalization $\|\En(c_h)-\bar{\En}(c_h)\|_{\infty,\Omega}\leq 1/2$. Therefore we can simplify the entropy viscosity \eqref{eq:eg:nu-e} to
\begin{align}
\label{eq:eg:simp-u}
\tilde\mu_\En &= 2\,\lambda_\En\,h_E^2\,\|R_h\|_{\infty,E} \leq \mu_\En.
\end{align}

We introduce an upper bound
\begin{align}
\label{eq:eg:nu-max}
\mu_{\max} = \lambda_{\max}\,h_E\,\|{f}'(c_h)\|_{\infty,E},
\end{align}
for some tunable parameter $\lambda_{\max}$, and set the entropy viscosity as
\begin{align}
\label{eq:eg:nu}
\mu_h = \min(\mu_{\max}, \tilde\mu_\En).
\end{align}
The original problem is modifed to
\begin{align}
\label{eq:eg:mod-cons}
\partial_t c + \div (F(c)) - \div(\mu_h\grad c) = 0.
\end{align}
Now for our equation \eqref{eq:2p:new-sat}, we adapt the ideas as follows.  We define the fully discrete entropy residual
\begin{align}\label{eq:ent-res-sw}
R_h^m &= \frac{(S_{w,h}^{m})^2 - (S_{w,h}^{m-1})^2}{2\delta t^m}
 + \phi^{-1}S_{w,h}^{m}f_w(S_{w,h}^{m})\u_a^m\cdot\grad S_{w,h}^{m}\\\nonumber
&\quad + \phi^{-1}S_{w,h}^{m}f_w(S_{w,h}^{m})\div\u_a^m,
\end{align}
locally on an element $E\in\cT_h$ (so the term $\grad S_{w,h}^{m}$ is well-defined).  In fact,
we drop the final term, since $\div\u_a^m=0$ except at wells. We define $\tilde\mu_h^m$ using \eqref{eq:eg:simp-u}, evaluating the maximum of $|R_h^m|$ by considering only the quadrature points used in the numerical integration formulas.  Similarly, we evaluate \eqref{eq:eg:nu-max}, which is now
\begin{align}
\label{eq:eg:nu-max-sw}
\mu_{\max}^m = \lambda_{\max}\,h_E\,\|\phi^{-1}S_{w,h}^{m} f_w'(S_{w,h}^{m})\u_a\|_{\infty,E},
\end{align}
and then we can set $\mu_h^m$ as in \eqref{eq:eg:nu}.

\section{IMPES coupling of flow and transport}\label{sec:impes}

We summarize the solution procedure, which is of the implicit pressure, explicit saturation (IMPES) type \cite{CHM_2006}. The wetting phase saturation $S_{w,h}^{0}\in X_h(\cT_h)$ is given as an approximation of the initial condition $S_0$ in \eqref{eq:2p:new-ic}, perhaps as the $L^2$-projection as stated in \eqref{eq:2p:trans1-ic}.  We set $S_{w,h}^{-1}=S_{w,h}^0$ to start the algorithm.

Now for time step $(t^m,t^{m+1}]$, $S_{w,h}^{m-1}$ and $S_{w,h}^{m}$ are given, and we do the following substeps.
\begin{enumerate}
\item Construct $\u_{c,h}^{m}\in\V_h$ and $\hat\zeta_h\in M_h$ by solving \eqref{eq:2p:uc-re1} and \eqref{eq:cap-flux-weak-continuity} using $S_{w,h}^{m}$ in place of $S_w$.  That is, solve \eqref{eq:2p:hf312} for the DoFs of $\hat\zeta_h$ and then define the DoFs of $\u_{c,h}^{m}$ using \eqref{eq:2p:uc-rec2}.

\item Solve the flow system \eqref{eq:2p:flow1}--\eqref{eq:2p:flow3} for $\u_{a,h}^m\in\V_h$, $\Phi_{w,h}^m\in W_h$, and $\hat\Phi_{w,h}^m\in M_h$ using the given saturation $S_{w,h}^{m}$ and the capillary flux $\u_{c,h}^m$ from the previous substep.  That is, solve for the DoFs using \eqref{eq:2p:hf17}--\eqref{eq:2p:hf25}.

\item For each element $E\in\cT_h$, compute the entropy viscosity $\mu_h^m$ from \eqref{eq:ent-res-sw},
\eqref{eq:eg:simp-u}, \eqref{eq:eg:nu-max-sw}, and  \eqref{eq:eg:nu}.  This requires $S_{w,h}^{m-1}$, $S_{w,h}^{m}$, and $\u_a^m$ from the previous substep.

\item Solve the stabilized transport equation \eqref{eq:2p:trans1} for $S_{w,h}^{m+1}\in X_h(\cT_h)$ using $S_{w,h}^{m}$ and, from previous substeps, $\mu_h^m$ and $\u_{a,h}^{m}$.  This can be done efficiently using the linear solution techniques developed in \cite{Lee_Lee_Wheeler_2016_eg}.
\end{enumerate}

\section{Extension to three dimensions}\label{sec:3d}

Our results extend to three dimensions, that is, to meshes of cuboidal hexahedra.  New families of mixed finite elements that approximate the velocity and pressure accurately were defined by the current authors in \cite{Arbogast_Tao_2018_atSpaces}, and they have the minimal number of DoFs for general hexahedra.  One could also use the macro-elements of Cockburn and Fu \cite{Cockburn_Fu_2017_mDecompIII}.  Unfortunately, at this writing, it is not known how to define serendipity elements that maintain accuracy on general hexahedra.  However, direct serendipity elements were defined in the Ph.D.\ dissertation of Tao \cite{Tao_2017_phd} for special hexahedra, namely, truncated cubes.  A truncated cube is a cuboidal hexahedron with two pairs of parallel faces (so the cross-section is parallelogram).

The dimension of the direct serendipity element of index $r$ on a truncated cube is $\dim\Po_r+3(r+1)$. Table~\ref{tab:eg:compare-3d} shows the number of DoFs for an $n\times n\times n$ mesh. Compared with using fully discontinuous spaces or standard mapped tensor product elements, the reduction in the number of degrees of freedom is much more significant than in two dimensions.  The enhanced Galerkin method is also dramatically better on hexahedra than using tetrahedra, since it takes at least five tetrahedra to fill a fixed hexahedron.

\begin{table}[ht]
\caption{A comparison of the global number of DoFs of
DG-$\Qu_r$, DG-$\Po_r$, EG-$\Qu_r$, and EG-$\DS_r$ on an $n\times n\times n$ mesh.}
\label{tab:eg:compare-3d}
{\centering\tabcolsep=4pt
\begin{tabular}{c|cccc}
\hline
$\myStrut r$ & DG-$\Qu_r$ & DG-$\Po_r$ & EG-$\Qu_r$ & EG-$\DS_r$\\
\hline
$\myStrut 0$ & $n^3$ & $n^3$ & --- & ---\\
$1$ & $8n^3$ & $4n^3$ & $2n^3+3n^3+3n+1$ & $2n^3+3n^2+3n+1$ \\
$2$ & $27n^3$ & $10n^3$ & $9n^3+12n^2+6n+1$ & $5n^3+9n^2+6n+1$ \\
$3$ & $64n^3$ & $20n^3$ & $28n^3+27n^2+9n+1$ & $8n^3+15n^2+9n+1$ \\
\hline 
\end{tabular}}
\end{table}

\section{Numerical tests}\label{sec:numericalResults}

In this work, the absolute permeability tensor is chosen as $\mathbf{K}=k\mathbf{I}$, where $\mathbf{I}$ is the identity matrix and $k>0$ is possibly heterogeneous over the domain.  We consider three numerical test cases.
They all use $s_{\text{form}}=-1$, i.e., the NIPG variant of the method.  We implemented the numerical methods within the deal.II framework \cite{BDHHKKMTW_2016_dealII84}.

\subsection{Coupled flow and transport}

We begin with a simple test of the numerical procedure, taken from \cite{Sun_Liu_2009_eg}, except that we use quadrilateral or hexahedral mesh elements instead of triangles.

Let $\lambda_t=1$, $f_w(S_w)=S_w$, and ignore capillary and gravity forces. The model reduces to a coupled linear flow and transport model that is used to represent single phase flow of the concentration $c=S_w$ of a dilute solute.  That is, we have
\begin{alignat}2
\label{eq:flow_darcy_c}
&\u = -\mathbf{K}\grad p &&\quad\text{in }\Omega,\\
\label{eq:flow_conservation_c}
&\div\u = q &&\quad\text{in }\Omega,\\
\label{eq:transport_c}
&\phi\frac{\partial c}{\partial t} + \div(c\u - \mu_h\grad c) = q_c&&\quad\text{in }\Omega.
\end{alignat}
We take $\phi = 1$, set $k=1$, and assume that there are no injection nor extraction sources within the domain (i.e., $q=q_c=0$).

\subsubsection{Two dimensions}

The problem is solved in a unit square with a circular hole in the middle. The flow equation \eqref{eq:flow_darcy_c}--\eqref{eq:flow_conservation_c} is given a Dirichlet boundary condition $p=1.5$ at the left side and $p=0.5$ at the right, and the top and bottom sides, as well as the circular hole, are given a no flow Neumann boundary condition ($\u\cdot\n=0$). The transport equation \eqref{eq:transport_c} is given an inflow boundary condition $c=1$ on the right side (i.e., on $\Gamma_{in}$), and the initial condition is $c=0$.

\begin{figure}
    \centering
        \includegraphics[width=.3\textwidth]{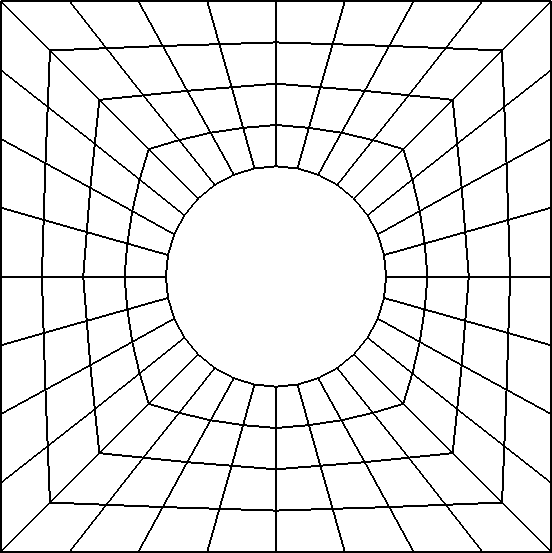}
    \caption{An example of the square domain with a hole. There are 128 convex quadrilateral
    elements in this coarse mesh.}\label{fig:eg:2d-cctn-grid2}
\end{figure}

Figure~\ref{fig:eg:2d-cctn-grid2} is an example of the mesh at refinement level $2$ with $128$ convex quadrilateral elements.  The actual simulation is conducted with a similar mesh at refinement level $5$ with $8192$ convex quadrilateral elements.

\begin{figure}
    \centering
    \begin{subfigure}[t]{0.36\textwidth}
        \centering
        \includegraphics[height=.79\textwidth]{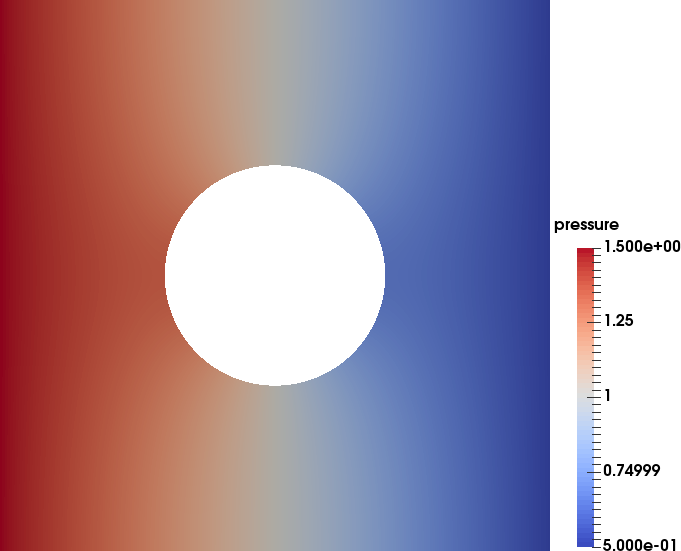}\ \ \quad
    \end{subfigure}%
    \qquad\qquad\qquad~
    \begin{subfigure}[t]{0.36\textwidth}
        \centering
        \includegraphics[height=.8\textwidth]{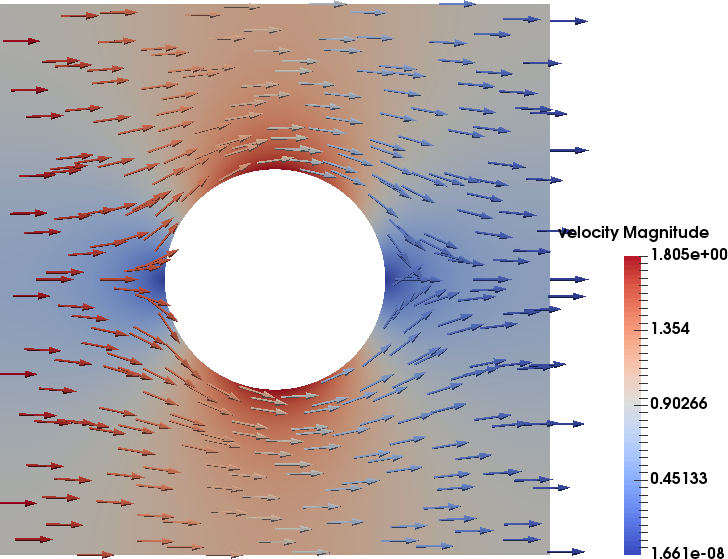}
    \end{subfigure}
    \caption{The pressure (top) and velocity field (bottom) of the coupled
    flow and transport system. The color indicates the magnitude of the velocity field,
    and the arrows indicate the direction.}
    \label{fig:eg:2d-cctn-pres-vel}
\end{figure}

The velocity $\u$ is approximated using the first order full $H(\Div)$-approximation Arbogast-Correa space ($\u_h\in\text{AC}_1$) \cite{Arbogast_Correa_2016}, which means that the pressure is approximated by a discontinuous polynomial of degree~1 (i.e., $p_h\in\bar X^P_1$). In Figure~\ref{fig:eg:2d-cctn-pres-vel}, we show the pressure and velocity for the flow problem \eqref{eq:flow_darcy_c}--\eqref{eq:flow_conservation_c}.

\begin{figure*}
    \centering
    \begin{subfigure}[t]{0.34\textwidth}
        \centering
        \includegraphics[width=\textwidth]{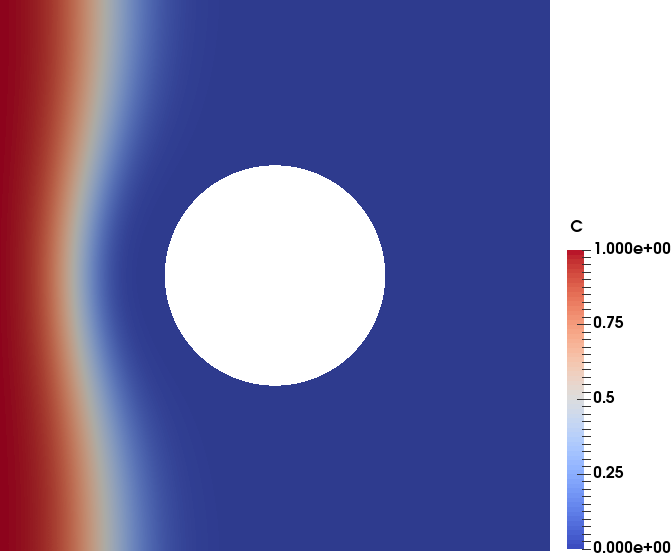}
        \vspace*{-16pt}\subcaption*{\kern-3em $m=200$}\vspace*{6pt}
    \end{subfigure}%
    \qquad\quad
    \begin{subfigure}[t]{0.34\textwidth}
        \centering
        \includegraphics[width=\textwidth]{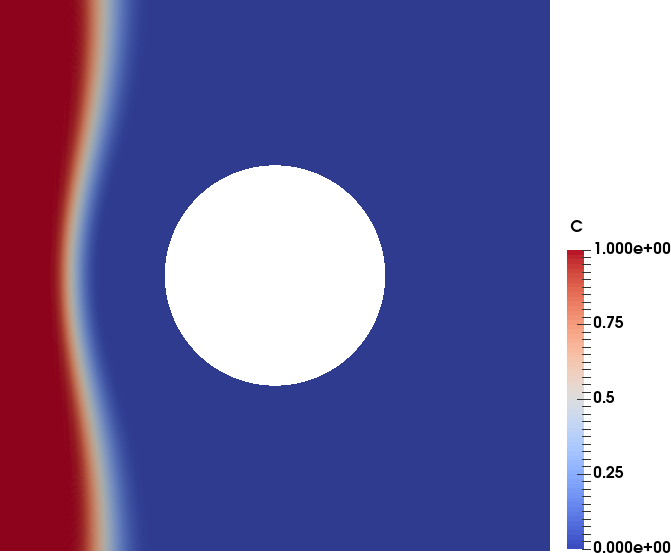}
        \vspace*{-16pt}\subcaption*{\kern-3em $m=200$}\vspace*{6pt}
    \end{subfigure}
    ~
    \begin{subfigure}[t]{0.34\textwidth}
        \centering
        \includegraphics[width=\textwidth]{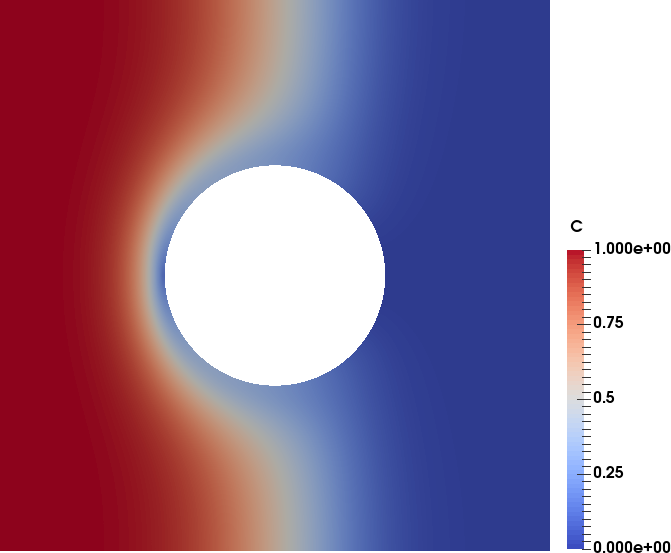}
        \vspace*{-16pt}\subcaption*{\kern-3em $m=500$}\vspace*{6pt}
    \end{subfigure}%
    \qquad\quad
    \begin{subfigure}[t]{0.34\textwidth}
        \centering
        \includegraphics[width=\textwidth]{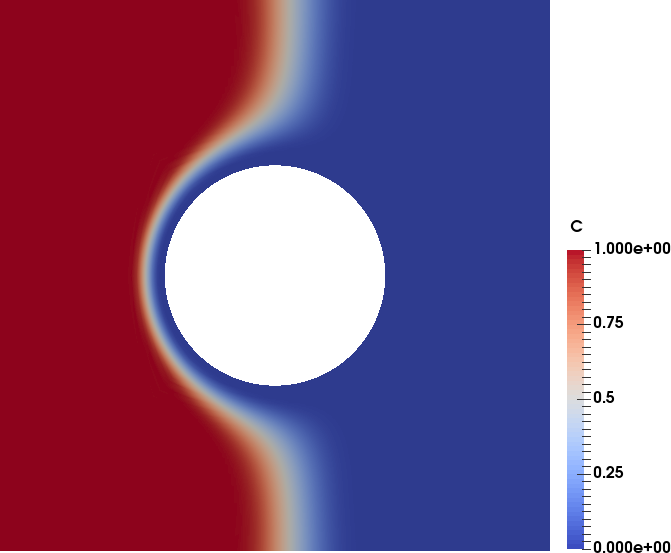}
        \vspace*{-16pt}\subcaption*{\kern-3em $m=500$}\vspace*{6pt}
    \end{subfigure}
    ~
        \begin{subfigure}[t]{0.34\textwidth}
        \centering
        \includegraphics[width=\textwidth]{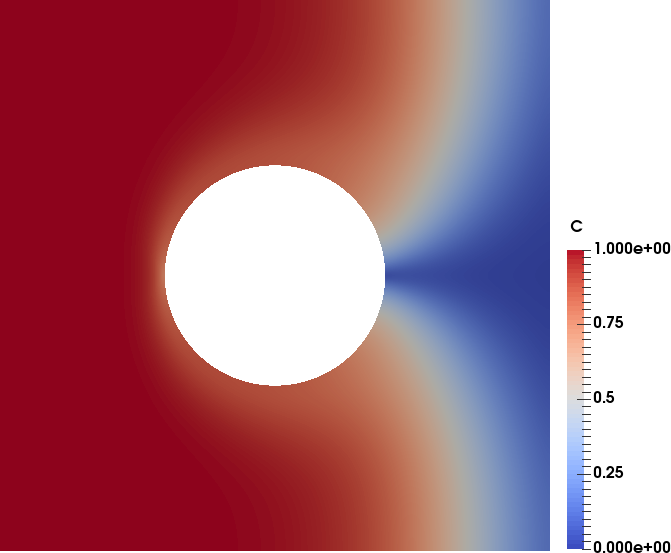}
        \vspace*{-16pt}\subcaption*{\kern-3em $m=800$}\vspace*{6pt}
    \end{subfigure}%
    \qquad\quad
    \begin{subfigure}[t]{0.34\textwidth}
        \centering
        \includegraphics[width=\textwidth]{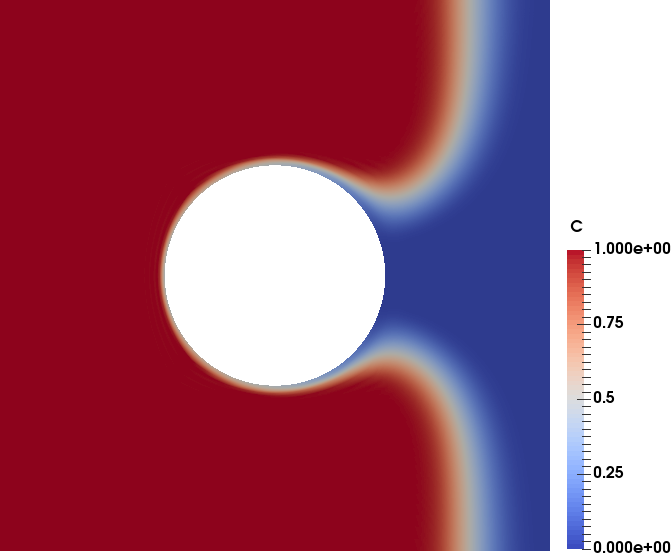}
        \vspace*{-16pt}\subcaption*{\kern-3em $m=800$}\vspace*{6pt}
    \end{subfigure}
    ~
    \begin{subfigure}[t]{0.34\textwidth}
        \centering
        \includegraphics[width=\textwidth]{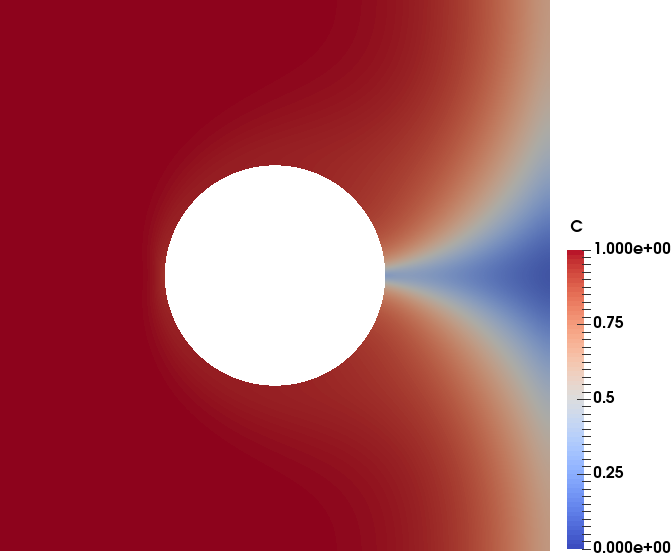}
        \vspace*{-16pt}\subcaption*{\kern-3em $m=1000$}\vspace*{6pt}
    \end{subfigure}%
    \qquad\quad
    \begin{subfigure}[t]{0.34\textwidth}
        \centering
        \includegraphics[width=\textwidth]{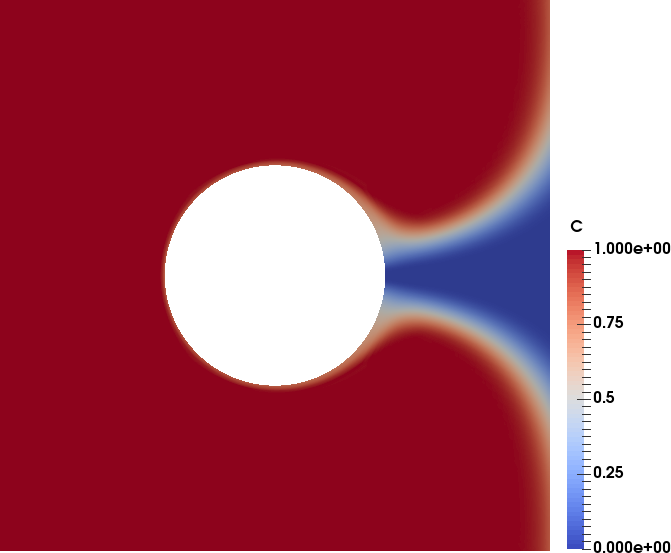}
        \vspace*{-16pt}\subcaption*{\kern-3em $m=1000$}\vspace*{6pt}
    \end{subfigure}
    \caption{Concentration values for each time step $m$. The left column is simulated with
    $\lambda_{\max}=1$, $\lambda_\En=\infty$ (so $\mu_h = \mu_{\max}$), and the right column is simulated with
    $\lambda_{\max}=1$, $\lambda_\En=0.5$ (so $\mu_h = \min(\mu_{\max}, \tilde\mu_\En)$).}
    \label{fig:eg:2d-cctn}
\end{figure*}    

The transport equation \eqref{eq:transport_c} is solved using the EG-$\DS_2$ method, i.e., we take $X_h(\cT_h) = X_2^{DS}(\cT_h)$ in~\eqref{eq:2p:trans1}. The time step is $\delta t = 0.1 h_{\min} =0.00105998$, where $h_{\min}$ is the minimal cell diameter of all elements in the mesh $\cT_h$.  The simplified version of entropy viscosity is applied.  The results are shown in Figure~\ref{fig:eg:2d-cctn}, where the left column show the concentration values at time step $m=200$, $500$, $800$, and $1000$ when only the maximum stabilization is applied ($\mu_h=\mu_{\max}$). The right column of Figure~\ref{fig:eg:2d-cctn} show the concentration values at the same time steps with the entropy stabilization ($\mu_h=\min(\mu_{\max}, \tilde\mu_E)$).  It is obvious that the concentration values with the entropy stabilization show a sharper front.

Our numerical results compare favorably in accuracy to those of Sun and Liu \cite{Sun_Liu_2009_eg}.  The advantage to our approach is that we use many fewer DoFs.  We achieve this by using quadrilateral elements and serendipity-based enriched Galerkin finite element spaces rather than triangular meshes.   This test was also performed in \cite{Lee_Lee_Wheeler_2016_eg}, and again we see comparable accuracy.  However, that work was constrained to use rectangular meshes for accuracy, and so the internal hole was modified to be a rectangle.  We were not subject to this restriction, since our methods are accurate on quadrilateral meshes and we can follow a curved interface with them (recall Figure~\ref{fig:eg:2d-cctn-grid2}).

\subsubsection{Three dimensions (3D)}

\begin{figure}
    \centering
        \includegraphics[width=.3\textwidth]{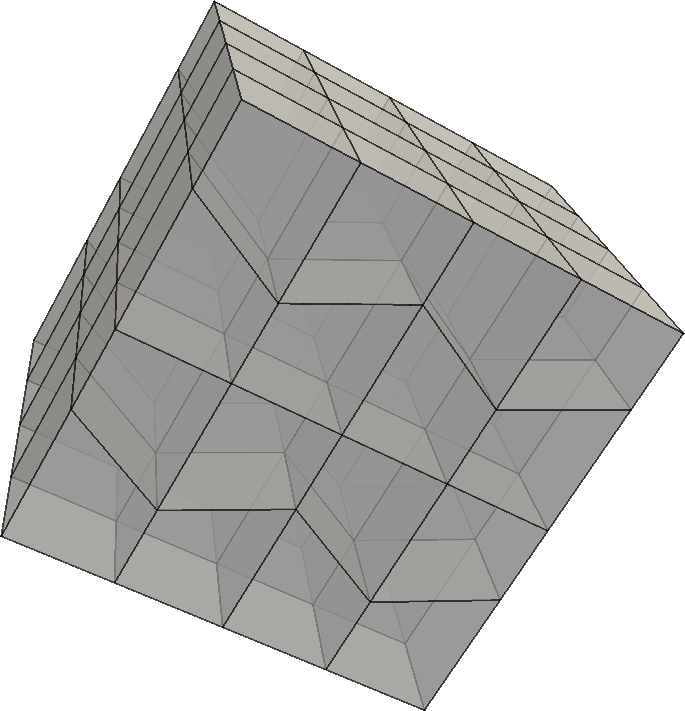}
    \caption{An example of the mesh for the 3D concentration flow at refinement level
    2 with $4\times 4\times 4$ truncated cubes.}\label{fig:eg:3d-cctn-grid2}
\end{figure}

\begin{figure}
    \centering
    \begin{subfigure}[t]{0.36\textwidth}
        \centering
        \includegraphics[height=.8\textwidth]{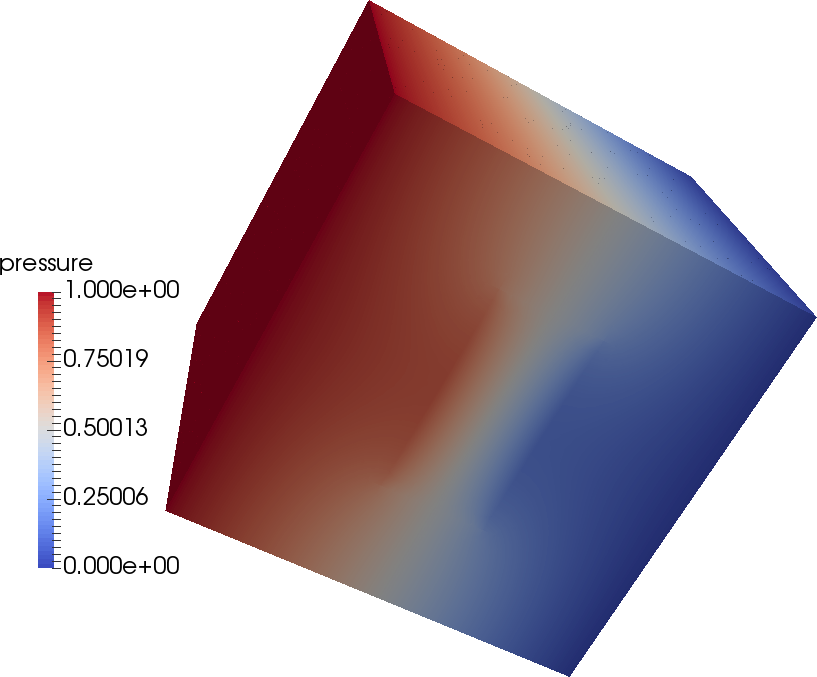}
    \end{subfigure}%
    \qquad\qquad\qquad~
    \begin{subfigure}[t]{0.36\textwidth}
        \centering
        \includegraphics[height=.8\textwidth]{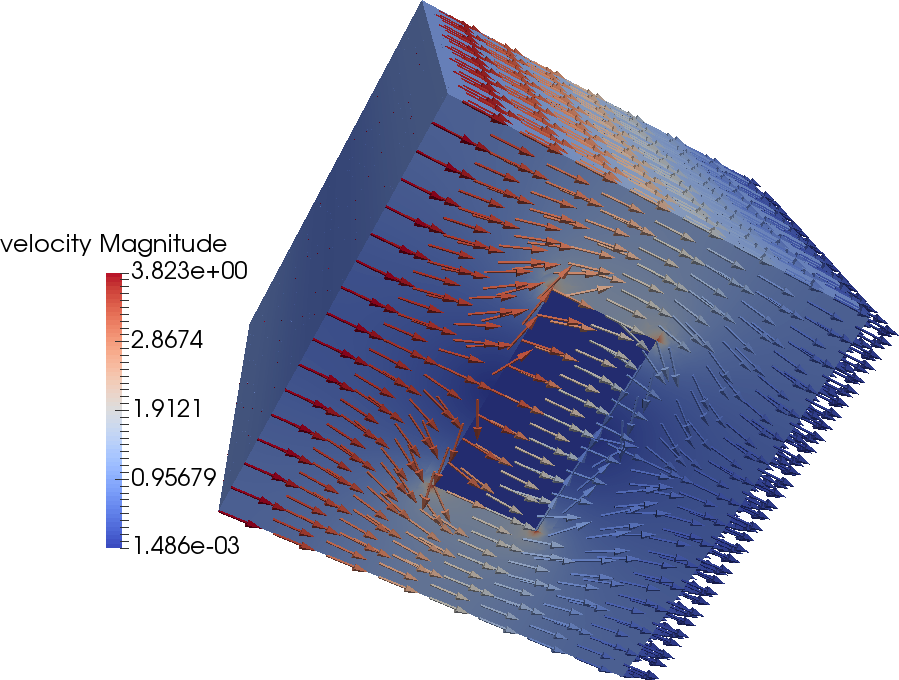}
    \end{subfigure}
    \caption{The pressure (top) and velocity field (bottom) of the coupled
    flow and transport system in 3D. The color indicates the magnitude
    of the velocity field and the arrows indicate the direction of flow.}
    \label{fig:eg:3d-cctn-pres-vel}
\end{figure}

\begin{figure*}
    \centering
    \begin{subfigure}[t]{0.32\textwidth}
        \centering
        \includegraphics[width=\textwidth]{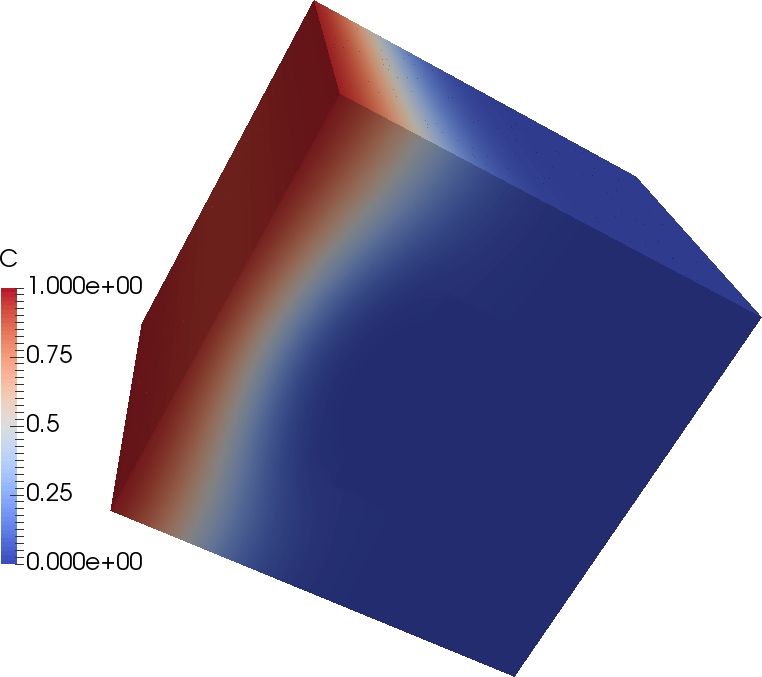}
        \vspace*{-24pt}\subcaption*{\kern-3em $m=40$}\vspace*{6pt}
    \end{subfigure}%
    \qquad\quad
    \begin{subfigure}[t]{0.32\textwidth}
        \centering
        \includegraphics[width=\textwidth]{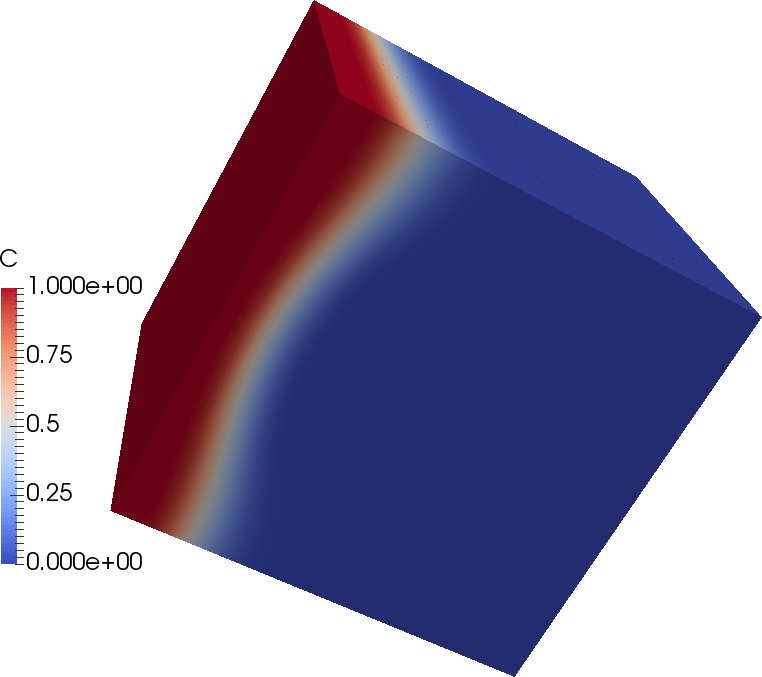}
        \vspace*{-24pt}\subcaption*{\kern-3em $m=40$}\vspace*{6pt}
    \end{subfigure}
    \medskip

    \begin{subfigure}[t]{0.32\textwidth}
        \centering
        \includegraphics[width=\textwidth]{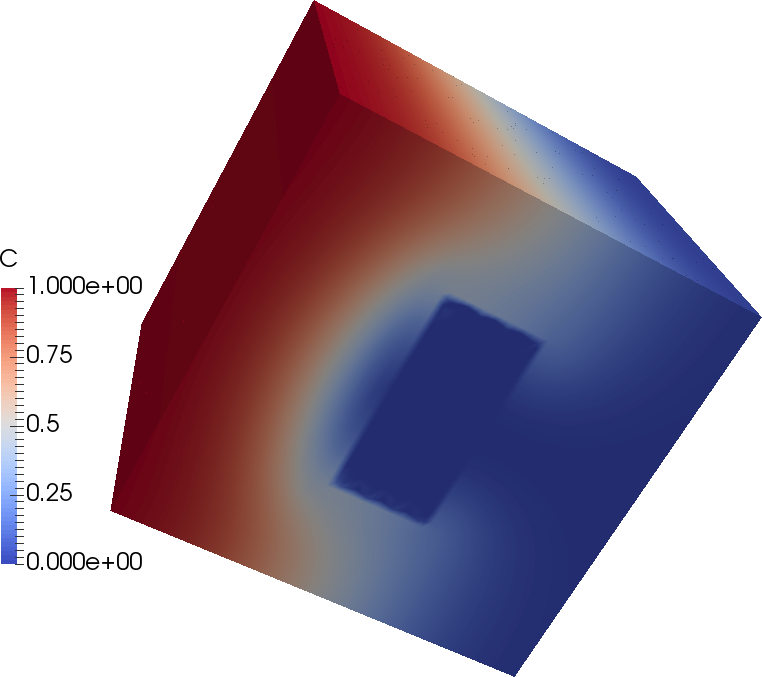}
        \vspace*{-24pt}\subcaption*{\kern-3em $m=100$}\vspace*{6pt}
    \end{subfigure}%
    \qquad\quad
    \begin{subfigure}[t]{0.32\textwidth}
        \centering
        \includegraphics[width=\textwidth]{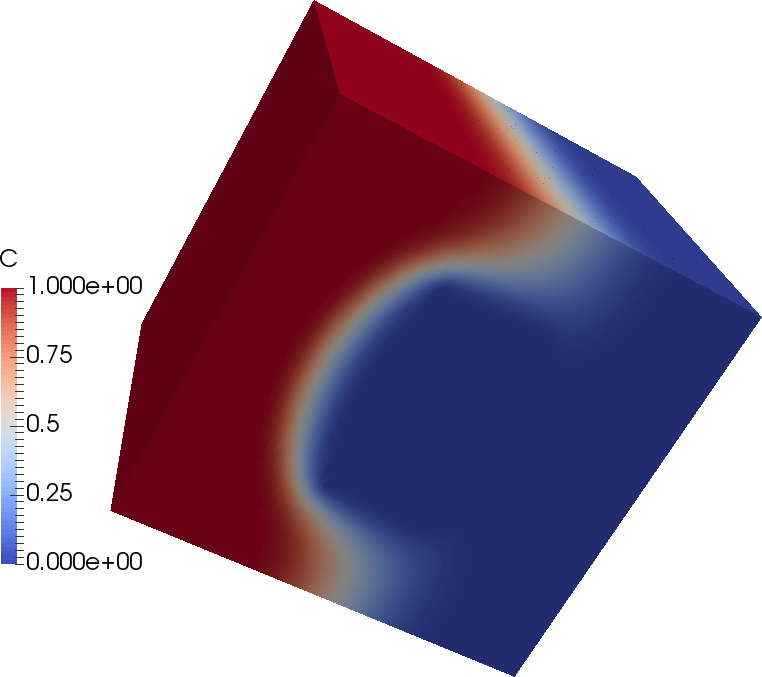}
        \vspace*{-24pt}\subcaption*{\kern-3em $m=100$}\vspace*{6pt}
    \end{subfigure}
    \medskip

        \begin{subfigure}[t]{0.32\textwidth}
        \centering
        \includegraphics[width=\textwidth]{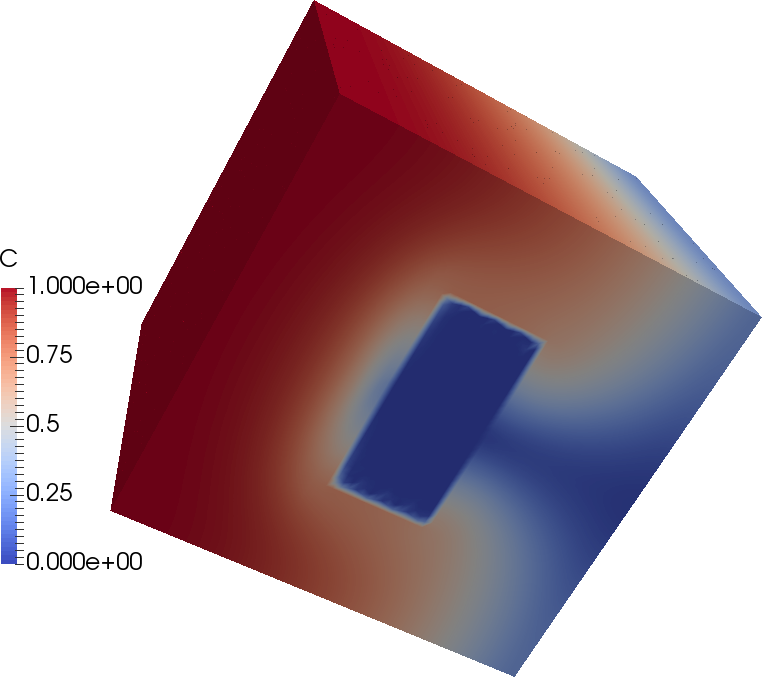}
        \vspace*{-24pt}\subcaption*{\kern-3em $m=160$}\vspace*{6pt}
    \end{subfigure}%
    \qquad\quad
    \begin{subfigure}[t]{0.32\textwidth}
        \centering
        \includegraphics[width=\textwidth]{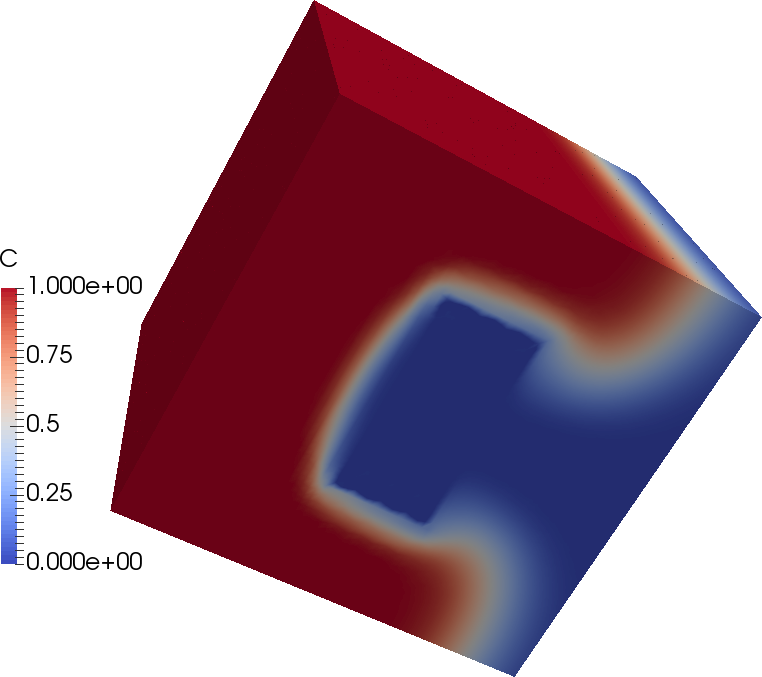}
        \vspace*{-24pt}\subcaption*{\kern-3em $m=160$}\vspace*{6pt}
    \end{subfigure}
    \medskip

    \begin{subfigure}[t]{0.32\textwidth}
        \centering
        \includegraphics[width=\textwidth]{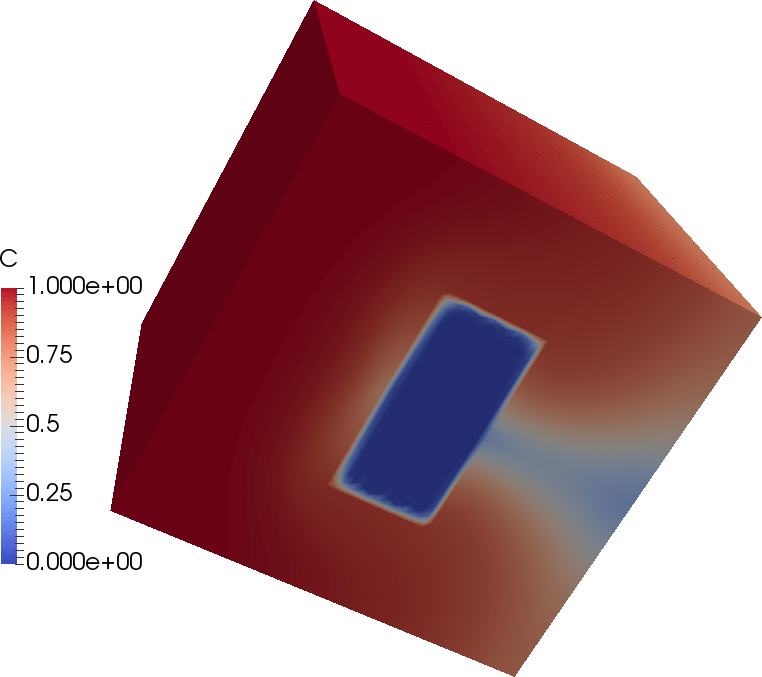}
        \vspace*{-24pt}\subcaption*{\kern-3em $m=250$}\vspace*{6pt}
    \end{subfigure}%
    \qquad\quad
    \begin{subfigure}[t]{0.32\textwidth}
        \centering
        \includegraphics[width=\textwidth]{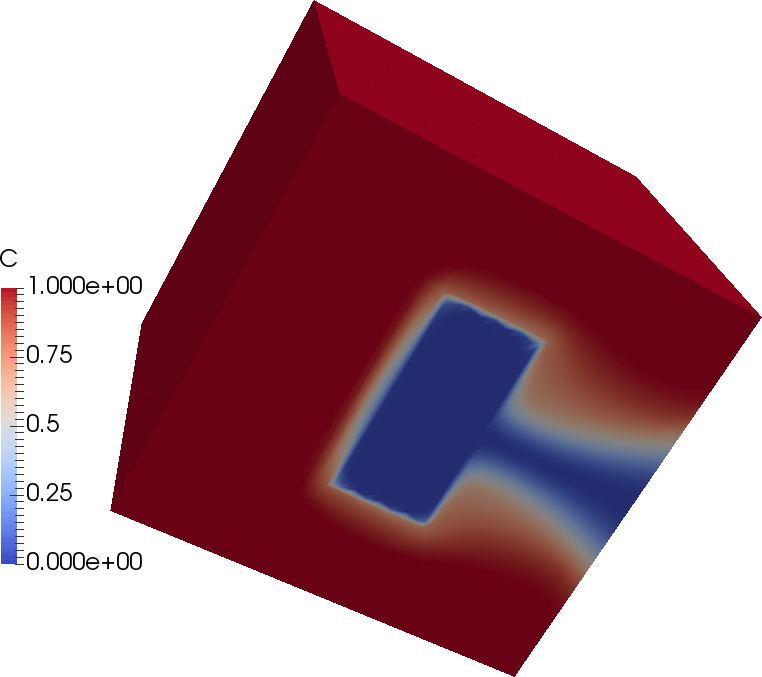}
        \vspace*{-24pt}\subcaption*{\kern-3em $m=250$}\vspace*{6pt}
    \end{subfigure}
    \caption{Concentration values for each time step $m$. The left column is simulated with     $\lambda_{\max}=1.0$, $\lambda_E=\infty$ (so $\mu_h = \mu_{\max})$, and the right column is simulated with $\lambda_{\max}=1.0$, $\lambda_E=0.5$ (so $\mu_h = \min(\mu_{\max}, \tilde\mu_\En)$).}
    \label{fig:eg:3d-cctn}
\end{figure*}   

For a 3D numerical simulation of \eqref{eq:flow_darcy_c}--\eqref{eq:transport_c}, the problem domain is $\Omega = [0,1]^3$. We put a low permeability block $\Omega_c = \left[3/8,5/8\right]\times [0,1]\times\left[1/4,3/4\right]$ in the middle where $k = 0.001$, and $k = 1$ in the rest of $\Omega$.  Figure~\ref{fig:eg:3d-cctn-grid2} is an example mesh at refinement level 2 with 64 truncated cubes, and the simulation is conducted on a similar triangulation at refinement level 5 with 32768 truncated cubes. The flow equation has a Dirichlet boundary condition $p=1.0$ at $(x,y,z)\in\{0.0\}\times[0,1]^2$ and $p=0.0$ at $(x,y,z)\in\{1.0\}\times[0,1]^2$; the other faces are set with a no flow Neumann boundary condition. There are no injection and extraction sources in the domain. The velocity field $\u$ is approximated in the first order full approximation Arbogast-Tao space (AT$_1$) defined in~\cite{Arbogast_Tao_2018_atSpaces} and the pressure is approximated by a discontinuous linear polynomial.  In Figure~\ref{fig:eg:3d-cctn-pres-vel}, we show the pressure and velocity for the flow problem.

For the transport equation, an entropy stabilized EG-$\DS_3$ method is applied, i.e., we take $X_h(\cT_h)=X_3^{DS}(\cT_h)$ in~\eqref{eq:2p:trans1}.  The transport equation is set with an inflow boundary condition $c=1$ on $\Gamma_{in}$ and the initial condition $c=0$. The time step is set to be $\delta t = 0.1 h_{\min} = 0.00541266$, where $h_{\min}$ is the minimal cell diameter of all the elements in the triangulation.  The porosity $\phi = 1$. The simplified entropy stabilization is taken as~\eqref{eq:eg:simp-u}.  The results are shown in Figure~\ref{fig:eg:3d-cctn} for time step $m=40$, $100$, $160$, $250$.  The left column shows the concentration values with the maximun stabilization ($\mu_h=\mu_{\max}$), and the right column show the values with the entropy stabilization ($\mu_h = \min(\mu_{\max}, \tilde\mu_\En)$).  We can observe that the concentration values with the entropy stabilization show sharper transport fronts and less numerical diffusion.

\subsection{Example 1: Effect of capillary pressure}

In this example, we show the effect of capillary pressure in a heterogeneous medium. The test is modified from that given in \cite{Hoteit_Firoozabadi_2008} in that we use quadrilateral mesh elements instead of rectangles.

\begin{figure}
    \centering
    \includegraphics[width=.35\textwidth]{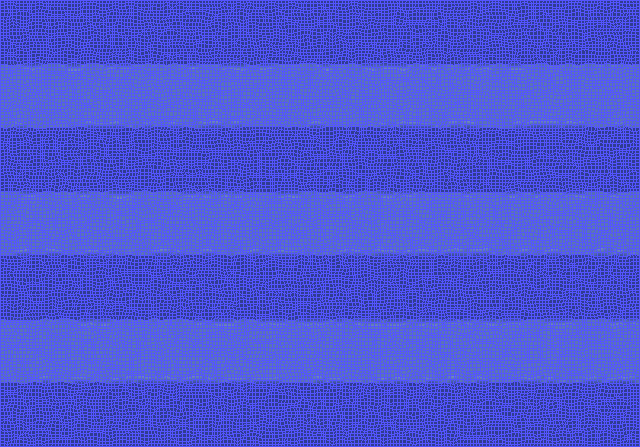}
    \caption{Computational domain for Example 1. There are $70\times 100$ 
    randomly distorted quadrilateral meshes and the permeability is 0.5
    in light zones and 0.01 in dark zones.}
    \label{fig:2p:layer-perm}
\end{figure}

This is a nondimensionalized model problem, so all physical units are ignored in the following description.  The computational domain is $\Omega=(0, 1.25)\times(0, 0.875)$. It is composed of layers of alternate permeabilities $k=0.5$ and $k=0.01$, as shown in Figure~\ref{fig:2p:layer-perm}, and rock porosity $\phi = 0.2$.  The fluid viscosities are $\mu_w = 1$ and $\mu_n = 0.45$.   The relative permeabilities are set as in~\eqref{eq:2p:krs} with $\beta=2$, the residual saturations are $S_{rw}=0$ and $S_{rn}=0$ (so $S_e = S_w$) and the capillary 
pressure is 
\begin{align}\label{eq:2p:pc001}
p_c(S_w) = \frac{-0.01}{\sqrt{k}}\log\big(\max(1.0, S_w+10^{-5})\big)
\end{align}
for the non-zero capillary pressure simulation, otherwise, it is set to be zero. We ignore gravity in this test, and therefore the densities are not required.

There are no injection and extraction sources inside the domain ($q_w=q_n=0$). For the flow equation, we set a Dirichlet boundary condition $\Phi_w = 1$ on the left hand side and $\Phi_w = 0$ on the right hand side, and a no flow Neumann boundary condition on the top and bottom sides of the domain. The inflow boundary for the saturation equation is the left hand side, and $S_B=1$ there. The domain is initially saturated with oil (the non-wetting phase), so $S_0 = 0$. 

The computational mesh $\cT_h$ consists of randomly distorted quadrilaterals, as shown in Figure~\ref{fig:2p:layer-perm}.  The finite element approximation spaces used are $\u_{a,h}\in\text{AC}_1$, $\Phi_{w,h}\in\bar X^P_1$, and $S_{w,h}\in X^{DS}_2$.  The time step is $\delta t=0.001$, and we use NIPG for the entropy stabilization ($\sigma_\gamma=1$). The entropy function $E(S_w) = \dfrac12 S_w^2$ and $F'(S_w) = S_w f_w'(S_w)\u_a$, and $\lambda_{\max} = 0.2$ and $\lambda_\En = 1$.

\begin{figure*}[p]
    \centering
    \bigskip
    \begin{subfigure}[t]{0.4\textwidth}
        \centering
        \includegraphics[width=\textwidth]{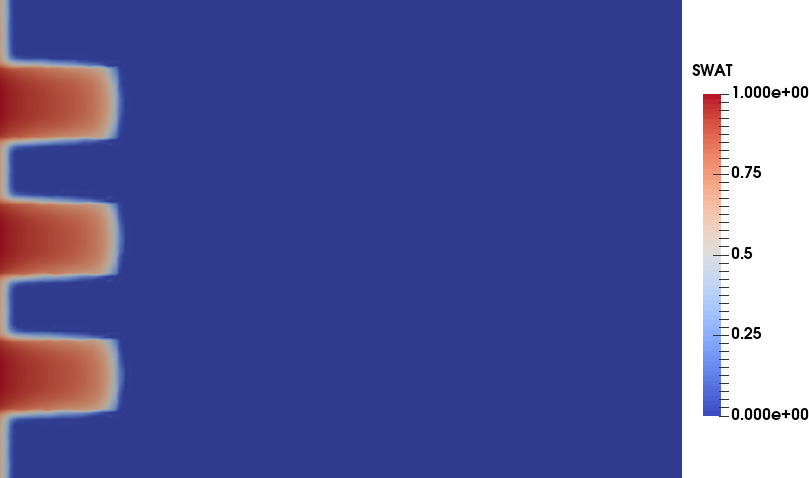}
        \vspace*{-16pt}\subcaption*{\kern-3em $m=50$}\vspace*{6pt}
    \end{subfigure}%
    \qquad
    \begin{subfigure}[t]{0.4\textwidth}
        \centering
        \includegraphics[width=\textwidth]{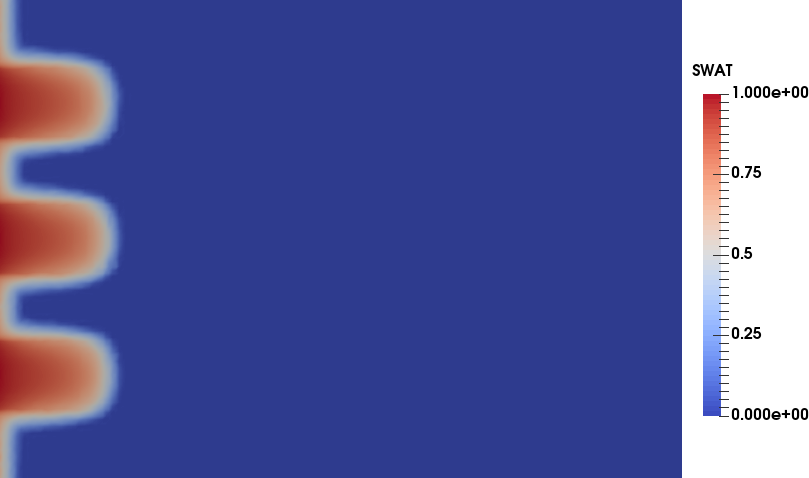}
        \vspace*{-16pt}\subcaption*{\kern-3em $m=50$ \\}\vspace*{6pt}
    \end{subfigure}
    \hfill
    \begin{subfigure}[h]{0.4\textwidth}
        \centering
        \includegraphics[width=\textwidth]{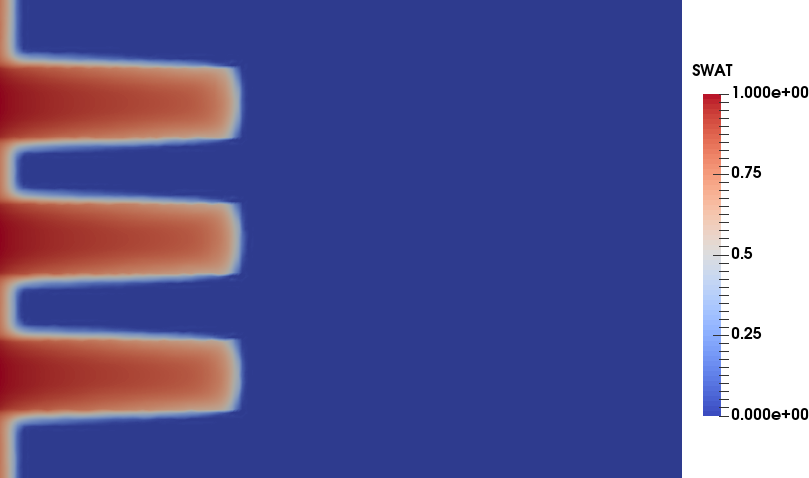}
        \vspace*{-16pt}\subcaption*{\kern-3em $m=125$}\vspace*{6pt}
    \end{subfigure}%
    \qquad
    \begin{subfigure}[h]{0.4\textwidth}
        \centering
        \includegraphics[width=\textwidth]{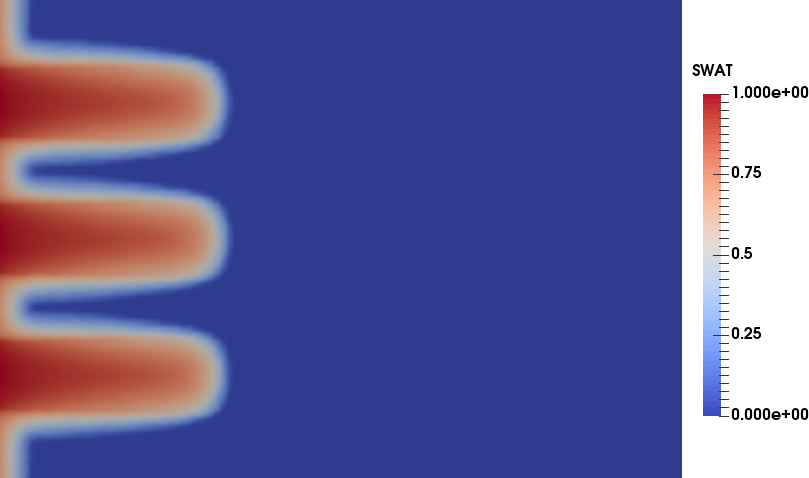}
        \vspace*{-16pt}\subcaption*{\kern-3em $m=125$ \\}\vspace*{6pt}
    \end{subfigure}
    \hfill
    \begin{subfigure}[h]{0.4\textwidth}
        \centering
        \includegraphics[width=\textwidth]{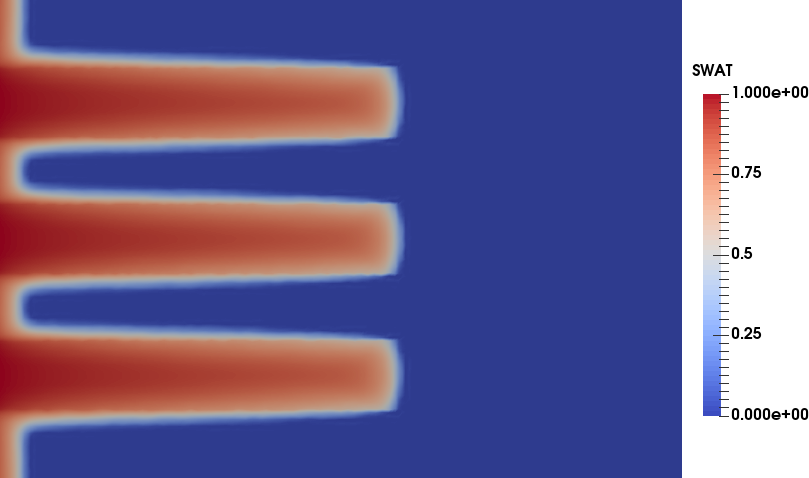}
        \vspace*{-16pt}\subcaption*{\kern-3em $m=250$}\vspace*{6pt}
    \end{subfigure}%
    \qquad
    \begin{subfigure}[h]{0.4\textwidth}
        \centering
        \includegraphics[width=\textwidth]{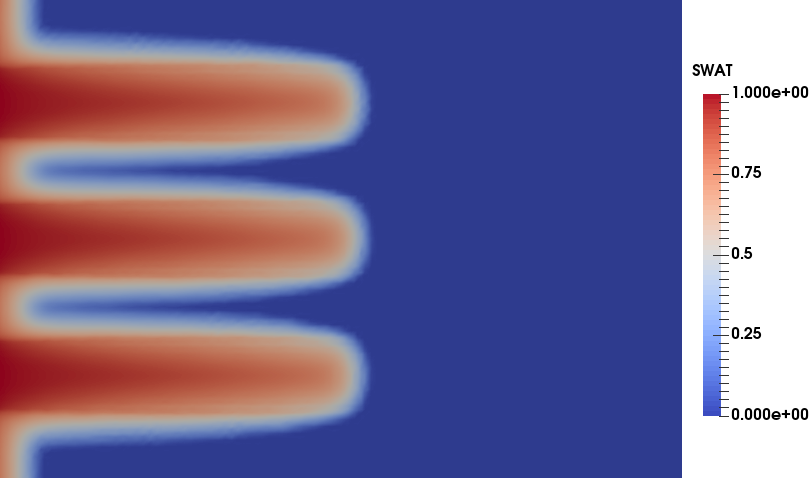}
        \vspace*{-16pt}\subcaption*{\kern-3em $m=250$ \\}\vspace*{6pt}
    \end{subfigure}
    \hfill
    \begin{subfigure}[h]{0.4\textwidth}
        \centering
        \includegraphics[width=\textwidth]{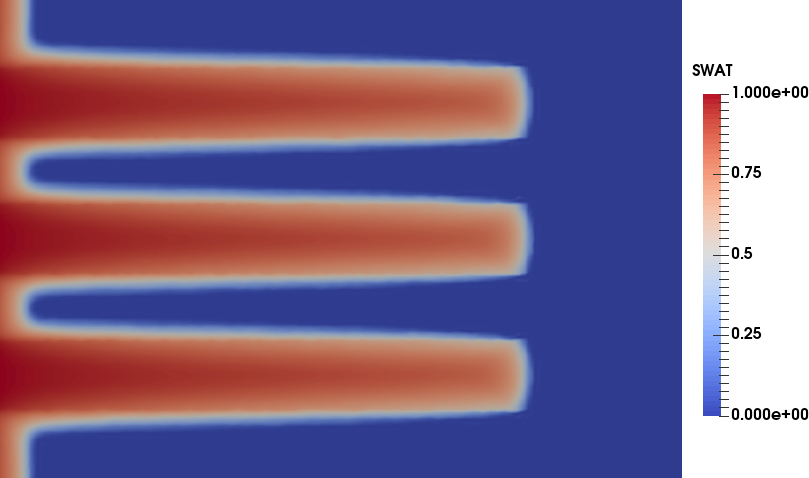}
        \vspace*{-16pt}\subcaption*{\kern-3em $m=375$}\vspace*{6pt}
    \end{subfigure}%
    \qquad
    \begin{subfigure}[h]{0.4\textwidth}
        \centering
        \includegraphics[width=\textwidth]{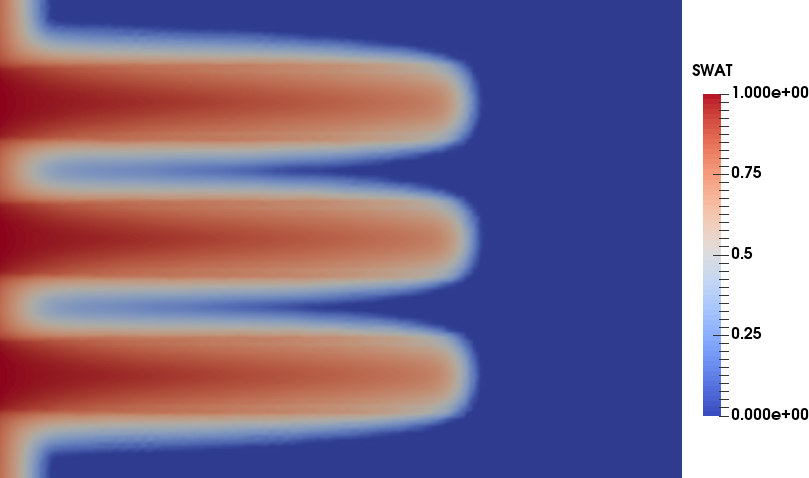}
        \vspace*{-16pt}\subcaption*{\kern-3em $m=375$}\vspace*{6pt}
    \end{subfigure}
    \hfill
    \caption{Wetting phase saturation profiles of Example 1 at time step $m=$ 50, 125, 250, and 375 with zero (left) and nonzero (right) capillary pressure.}
    \label{fig:2p:layer}
\end{figure*}

In Figure~\ref{fig:2p:layer}, we compare the simulated wetting phase saturation with and without capillary pressure at the same time steps.  We observe that the front of the saturation profile with capillary pressure slows down and is more diffuse compared to the zero capillary pressure case. Moreover, the velocity is much greater near the interface between different rock types due to the presence of heterogeneous capillarity, as similar numerical tests in~\cite{Lee_Wheeler_2018,Hoteit_Firoozabadi_2008,AJPW_2013,Kou_Sun_2010} show.

Our results compare favorably to those in \cite{Hoteit_Firoozabadi_2008}, where rectangular elements were used.  The test shows that grid distortion does not degrade the accuracy of the results.

\subsection{Example 2: Heterogeneous capillary pressure}

\begin{figure}[!htb]
    \centering
    \includegraphics[width=.4\textwidth]{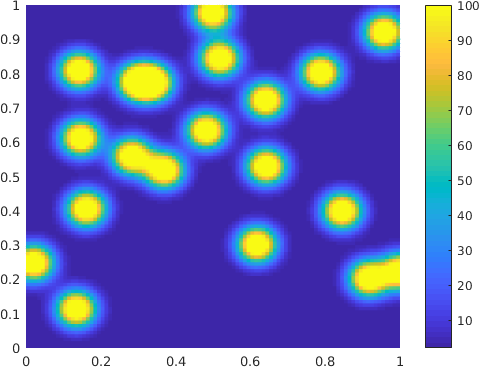}
    \caption{Computational domain for Example 2. 
    The color indicate the value of $k^{-1}$, where the permeability 
    $\mathbf{K}=k\mathbf{I}$. The minimal value of $k^{-1}$ is 2 and the maximal
    is 100.}
    \label{fig:2p:rand-perm}
\end{figure}

In this second example of two-phase flow, we show the difference of homogeneous and heterogeneous capillary pressure in a heterogeneous medium. Similar to the previous example, this is a nondimensionalized model problem so that all physical units are omitted. The computational domain is $\Omega=(0, 1)^2$.  A random permeability field is distributed over $\Omega$ as depicted in Figure~\ref{fig:2p:rand-perm}. The rock porosity $\phi$, the fluid viscosities $\mu_w$ and $\mu_n$, the relative permeabilities, the residual saturations, and initial saturation $S_0$ are the same as in Example~1, and gravity is ignored here as well. The heterogeneous capillary pressure is given in~\eqref{eq:2p:pc001}, and the homogeneous capillary pressure is set as
\begin{align}\label{eq:2p:pc002}
p_c(S_w) = -0.01\log\big(\max(1.0, S_w+10^{-5})\big),
\end{align}
which does not depend on the rock permeability.

There are no injection and extraction sources inside the domain ($q_w=q_n=0$).  The water (wetting phase) injection boundary is $\{0\}\times(0,0.05)$, where $\Phi_w = 1$, and the extraction boundary is $\{1\}\times(0.95,1)$, where $\Phi_w=0$. The rest of the flow boundary is set with a no flow Neumann boundary condition. The inflow boundary is the injection boundary, and $S_B=1$ there.

\begin{figure}[!htb]
    \centering
    \includegraphics[width=.25\textwidth]{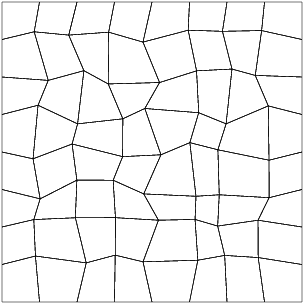}
    \caption{An example of a $8\times 8$ randomly distorted quadrilateral mesh.}
    \label{fig:rand-grid}
\end{figure}

The computational mesh has $100\times100$ distorted quadrilateral elements (as in the coarse mesh of Figure~\ref{fig:rand-grid}). The time step, approximation spaces, entropy stabilization and the NIPG formulation are the same as Example~1, except $\lambda_\En=1.5$ in this example.

\begin{figure*}
    \centering
    \bigskip
    \begin{subfigure}[t]{0.4\textwidth}
        \centering
        \includegraphics[width=\textwidth]{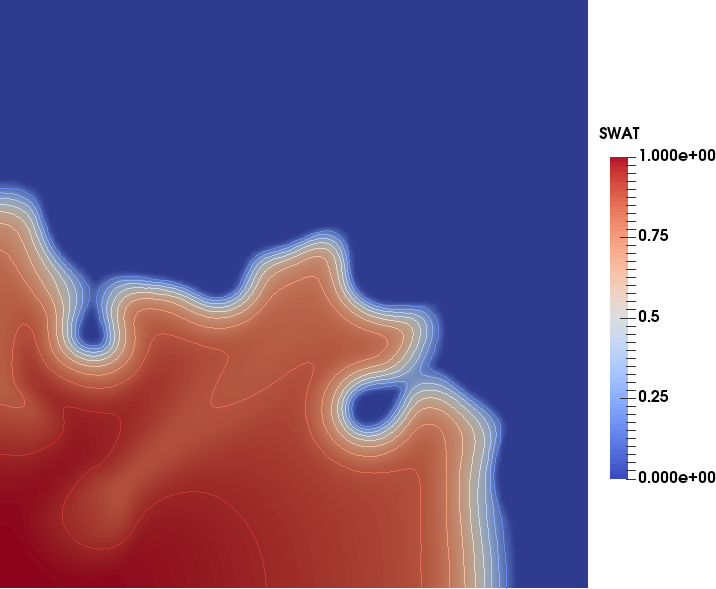}
        \vspace*{-16pt}\subcaption*{\kern-3em $m=1000$}\vspace*{6pt}
    \end{subfigure}%
    \qquad
    \begin{subfigure}[t]{0.4\textwidth}
        \centering
        \includegraphics[width=\textwidth]{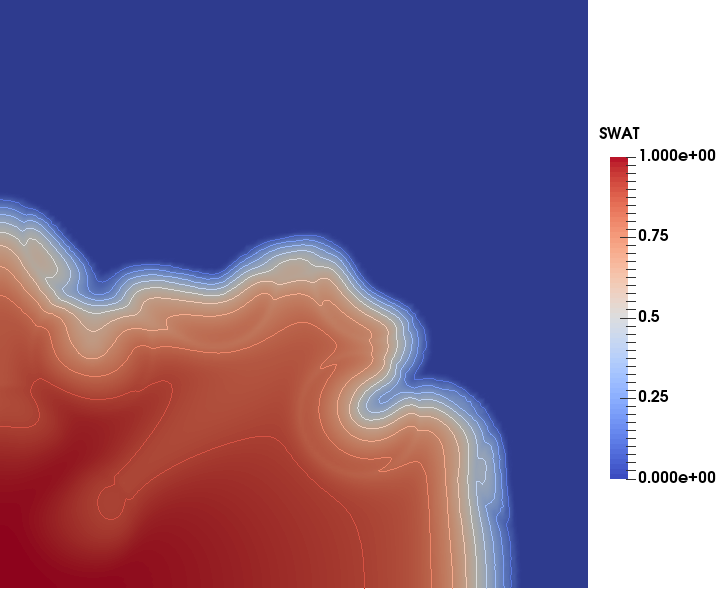}
        \vspace*{-16pt}\subcaption*{\kern-3em $m=1000$}\vspace*{6pt}
    \end{subfigure}
    ~
    \begin{subfigure}[t]{0.4\textwidth}
        \centering
        \includegraphics[width=\textwidth]{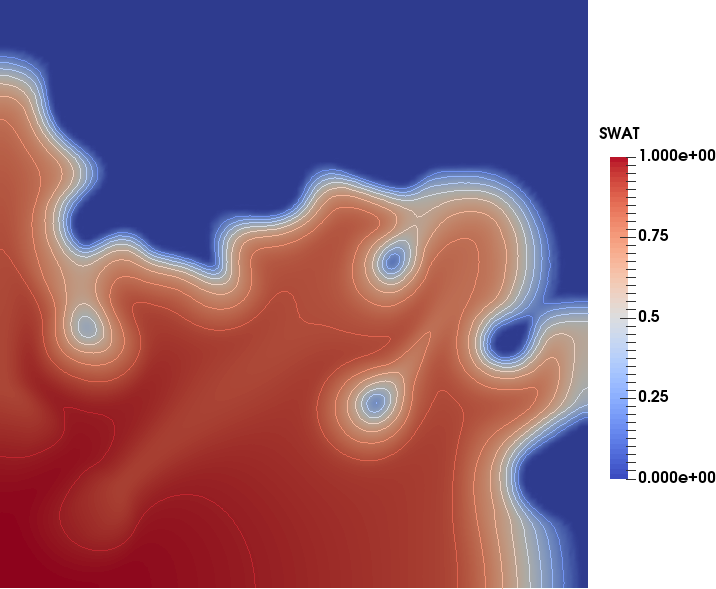}
        \vspace*{-16pt}\subcaption*{\kern-3em $m=1500$}\vspace*{6pt}
    \end{subfigure}%
    \qquad
    \begin{subfigure}[t]{0.4\textwidth}
        \centering
        \includegraphics[width=\textwidth]{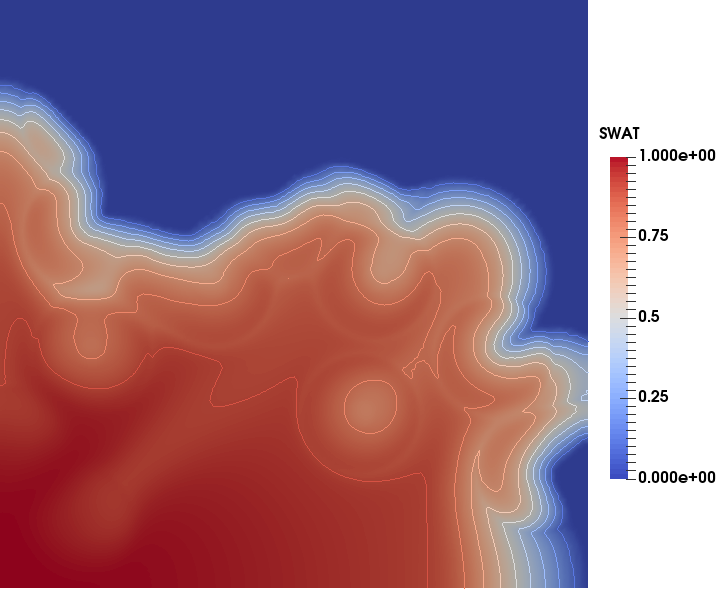}
        \vspace*{-16pt}\subcaption*{\kern-3em $m=1500$}\vspace*{6pt}
    \end{subfigure}
    ~
    \begin{subfigure}[t]{0.4\textwidth}
        \centering
        \includegraphics[width=\textwidth]{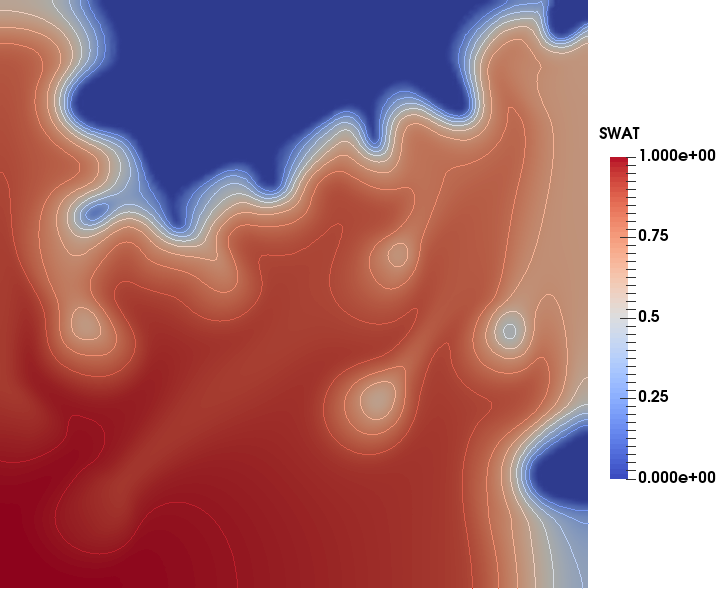}
        \vspace*{-16pt}\subcaption*{\kern-3em $m=2000$}\vspace*{6pt}
    \end{subfigure}%
    \qquad
    \begin{subfigure}[t]{0.4\textwidth}
        \centering
        \includegraphics[width=\textwidth]{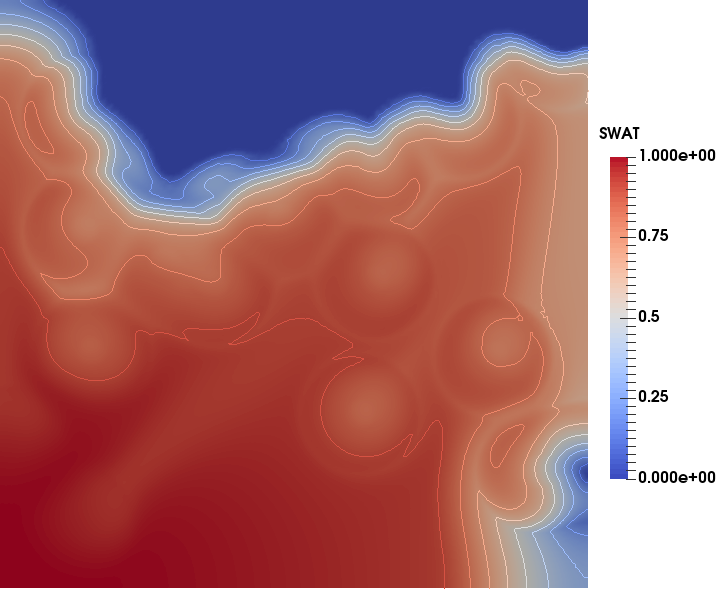}
        \vspace*{-16pt}\subcaption*{\kern-3em $m=2000$}\vspace*{6pt}
    \end{subfigure}
    \caption{Wetting phase saturation profiles at $m=$1000, 1500, and 2000 with homogeneous (left) and 
    heterogeneous (right) capillary pressure.
    There are 9 contour lines set at saturation levels $0.1$--$0.9$, $0.1$ step for each.}
    \label{fig:2p:rand}
\end{figure*}

Figure~\ref{fig:2p:rand} shows the profiles of the wetting phase saturation with homogeneous and heterogeneous capillary pressures at the same time steps.  At time step $m=1000$, there are two low permeability regions that are not invaded by the wetting phase if a homogeneous capillary pressure is applied.  On the other hand, with a heterogeneous capillary pressure, the wetting phase invades these regions.  At time step $m=2000$, the water front has almost reached the extraction boundary. The homogeneous capillary pressure case shows an unrecovered region near the bottom right of the domain, while that region is much better saturated if the heterogeneous capillary pressure is applied.

Being able to use quadrilateral meshes greatly aids the set-up of this problem.  It allows us to accurately model the shape of the low permeability inclusions.  We maintain logically rectangular mesh indexing, and we use many fewer DoFs than triangles would provide.

\section{Summary and conclusions}\label{sec:conc}

We described recently developed mixed and direct serendipity finite elements \cite{Arbogast_Correa_2016,Arbogast_Tao_2018x_serendipity}.  These elements provide accurate approximations on meshes of quadrilaterals and maintain the minimal number of degrees of freedom (DoFs) needed to maintain finite element conformity.  Elements that approximate to any order are available.  We developed a new enriched Galerkin method \cite{Sun_Liu_2009_eg} on quadrilaterals based on the direct serendipity finite elements.

We applied the new elements to computational simulation of two-phase flow in porous media.  We used the Hoteit-Firoozabadi formulation \cite{Hoteit_Firoozabadi_2008}, which separates the system in terms of the advective and capillary velocities $\u_a$ and $\u_c$ into an elliptic pressure equation for the flow velocity $\u_a$ and a hyperbolic saturation equation for $S_w$. A standard IMPES solution procedure allows the two parts of the system to be solved independently.

The pressure equation was discretized on quadrilateral meshes using a mixed finite element approximation and the accurate Arbogast-Correa spaces \cite{Arbogast_Correa_2016}.  However, the formulation requires construction of the divergence of the capillary flux.  We provided a novel implementation that does not break down when the system degenerates (i.e., one of the saturations tends to the residual value).

The saturation equation was also discretized accurately on quadrilateral meshes using the new enriched Galerkin method employing the direct serendipity spaces of the current authors \cite{Arbogast_Tao_2018x_serendipity}.  Being hyperbolic, the system needs to be stabilized, and we used the entropy stabilization procedure of Guermond, Pasquetti, and Popov \cite{Guermond_Pasquetti_Popov_2011_entropyVisc}, adapted to our specific equation.

Extension to three space dimensions is straightforward, up to the definition of the minimal DoF finite element spaces.  The mixed spaces have been defined in \cite{Arbogast_Tao_2018_atSpaces,Cockburn_Fu_2017_mDecompIII}, and the direct serendipity spaces were defined in \cite{Tao_2017_phd} for truncated cubes.

Numerical results showed that accurate results are obtained by our numerical method, in both two and three dimensions.  The ability to use quadrilateral and hexahedral meshes greatly increases one's ability to set up appropriate meshes for the problem to be solved. Our method uses a minimal number of DoFs, and so is quite efficient.  The low number of DoFs is due to two sources.  First, quadrilateral and hexahedral meshes use far fewer elements than ones based on triangles or simplices.  Second, we use minimal DoF finite element spaces, which greatly reduce the number of DoFs compared to, say ABF or Devloo et al.\ mixed finite element spaces \cite{ABF_2002,Bergot_Durufle_2013,Siqueira_Devloo_Gomes_2013} and standard enriched Galerkin methods based on tensor product polynomials.  Moreover, the method is efficient because, in many problems, we can maintain a logically rectangular mesh indexing.

Our numerical results also showed that our novel construction of the divergence of the capillary flux captures well the effects of capillary pressure.


\let\v=\vv


\end{document}